\documentclass[a4paper]{article}
\usepackage[margin=1in]{geometry}
\usepackage{amssymb}
 \usepackage{wrapfig}
 \usepackage{flafter}
 \usepackage{leftidx}
 \usepackage{mathrsfs}
 \usepackage{amssymb,amsmath}
 \usepackage{bbding}
\usepackage{fancyhdr}%调用排印页眉与页脚的宏包
\usepackage{mathrsfs}
\usepackage{stmaryrd}
\usepackage{amsfonts}
\usepackage{latexsym}
\usepackage{psfrag}
\usepackage{graphicx,subfigure}
\usepackage{dsfont}
\usepackage{bbm}
\usepackage{booktabs}
\usepackage{color,xcolor}
\usepackage{authblk}
\usepackage{appendix}
\usepackage{cite}
\allowdisplaybreaks

\usepackage{graphicx}
\usepackage{natbib,color}

\newcommand{\eps}{\varepsilon}
\newcommand{\nn}{\nonumber}
\newcommand{\bx}{{\bf x}}
\newcommand{\bT}{{\bf T}}
\newcommand{\bF}{{\bf F}}

\newtheorem{theorem}{Theorem}[section]
\newtheorem{lemma}[theorem]{Lemma}
\newtheorem{remark}[theorem]{Remark}
\graphicspath{{./}}

\begin{document}

\title{Improved uniform error bounds on time-splitting methods for the long-time dynamics of the weakly nonlinear Dirac equation}
% Short title for running heads:
%\shorttitle{Improved uniform error bounds on TSFP for NLDE}

\author{%
{\sc Weizhu Bao}\thanks{ Email: matbaowz@nus.edu.sg}
 \\
Department of Mathematics, National University of Singapore, Singapore 119076, Singapore\\
{\sc Yongyong Cai}\thanks{ Email: yongyong.cai@bnu.edu.cn}
 \\
Laboratory of Mathematics and Complex Systems and School of Mathematical Sciences, Beijing Normal University, Beijing 100875, China\\
{\sc and}\\
{\sc Yue Feng}\thanks{Email: fengyue@u.nus.edu}\\
Department of Mathematics, National University of Singapore, Singapore 119076, Singapore
}
% Short list of authors for running heads:
%\shortauthorlist{W. Bao \emph{et al.}}

\maketitle

\begin{abstract}
% Body of abstract:
{Improved uniform error bounds on time-splitting methods are rigorously proven for the long-time dynamics of the weakly nonlinear Dirac equation (NLDE), where the nonlinearity strength is characterized by a dimensionless parameter $\eps \in (0, 1]$ . We adopt a second order Strang splitting method to discretize the NLDE in time and combine the Fourier pseudospectral method in space for the full-discretization. By employing the {\sl regularity compensation oscillation} (RCO) technique where the high frequency modes are controlled by the regularity of the exact solution and the  low frequency modes are analyzed by phase cancellation and energy method, we establish improved uniform error bounds at $O(\eps^2\tau^2)$ and $O(h^{m-1}+ \eps^2\tau^2)$ for the second-order Strang splitting semi-discretizaion and full-discretization up to the long-time $T_{\eps} = T/\eps^2$ with $T>0$ fixed, respectively. Furthermore, the numerical scheme and error estimates are extended to an oscillatory NLDE which propagates waves with $O(\eps^2)$ wavelength in time and at $O(\eps^{-2})$ wave speed in space. Finally, numerical examples verifying our analytical results are given.}
% Keywords:
{nonlinear Dirac equation; long-time dynamics; time-splitting method; improved uniform error bound; regularity compensation oscillation (RCO).}
\end{abstract}

\section{Introduction}
\label{sec;introduction}
Long-time dynamics of Hamiltonian partial differential equation (PDE) has attracted much interest in recent years from both analytical and numerical aspects \citep{CV,FG,FGP,SH}. The change of a small nonlinear perturbation to the dynamics of a linear Hamiltonian PDE deserves careful considerations \citep{GL,HL}. The nonlinear Dirac equation (NLDE) is widely used in many fields such as the electron self-interaction \citep{Dirac1,ES}, quantum field theory \citep{TWE,Soler}, Bose-Einstein condensate (BEC) \citep{HC} and graphene as well as other 2D materials \citep{BHM,FW}. In this paper, we study the dynamics of the weakly NLDE in one or two dimensions (1D or 2D), which can be expressed for the two-component wave function $\Phi := \Phi(t, \bx) = (\phi_1(t, \bx), \phi_2(t, \bx))^T \in {\mathbb C}^2$ as \citep{Dirac1,Dirac2}
\begin{equation}
\label{eq:NLDE_21}
i\partial_t\Phi =  \Big(- i\sum_{j = 1}^{d}\sigma_j\partial_j + \sigma_3 \Big)\Phi +\eps^2 {\bF}(\Phi)\Phi, \quad x \in \Omega, \quad t > 0,
\end{equation}
with the initial data
\begin{equation}
\label{eq:initial}
\Phi(t=0, {\bf x}) = \Phi_0({\bf x}), \quad {\bf x} \in \overline{\Omega},
\end{equation}
where $\Omega = \prod^d_{j= 1} (a_j, b_j) \subset \mathbb{R}^d$ ($d = 1, 2$) is a bounded domain equipped with periodic boundary conditions. Here, $i = \sqrt{-1}$ is the imaginary unit, $t$ is time, $\mathbf{x}=(x_1,\ldots, x_d)^T$ is the spatial coordinate vector, $\partial_j$ represents $\partial_{x_j} $ for $j = 1, \ldots, d$, $\eps \in (0, 1]$ is a dimensionless parameter. $\sigma_1$, $\sigma_2$ and $\sigma_3$ are the Pauli matrices given as
	\begin{equation}
	\label{Pauli}
	\sigma_1 = \begin{pmatrix} 0 &\  1 \\ 1 &\   0 \end{pmatrix}, \quad
	\sigma_2 = \begin{pmatrix} 0 & \  -i \\ i &\    0\end{pmatrix}, \quad
	\sigma_3 = \begin{pmatrix} 1 &\  0 \\ 0 &\   -1 \end{pmatrix}.
	\end{equation}
The matrix nonlinearity ${\bf F}(\Phi) $ in \eqref{eq:NLDE_21} is usually taken as
 \begin{equation}
 {\bF}(\Phi) = \lambda_1(\Phi^{\ast}\sigma_3\Phi)\sigma_3 + \lambda_2|\Phi|^2 I_2, \quad \Phi \in \mathbb{C}^2,
 \label{eq:nonlinear}
 \end{equation}
where $\lambda_1$, $\lambda_2$ are two given real constants, $I_2$ is the $2 \times 2$ identity matrix and $|\Phi|^2 = \Phi^{\ast}\Phi$ with $\Phi^{\ast} = \bar{\Phi}^T$ representing the complex conjugate transpose of $\Phi$. The choice of the nonlinearity ${\bf F}(\Phi)$ is motivated by the so-called Soler model in quantum field theory, e.g., $\lambda_1 \neq 0$ and $\lambda_2 = 0$  \citep{FS,TWE,Soler}, and the BEC with a chiral confinement and/or spin-orbit coupling, e.g., $\lambda_1=0$ and $\lambda_2 \neq 0$ \citep{CEL,HC}.

The NLDE \eqref{eq:NLDE_21} is dispersive and time symmetric \citep{BCJY}. In addition, it  conserves the total mass as
\begin{equation}
\left\|\Phi(t, \cdot)\right\|^2	:= \int_{\Omega} |\Phi(t, \bx)|^2 d\bx = \int_{\Omega} \sum^2_{j = 1}|\phi_j(t, \bx)|^2 d\bx \equiv \left\|\Phi(0, \cdot)\right\|^2 = \left\|\Phi_0\right\|^2, \quad t \geq 0,
\end{equation}
and the energy as
\begin{equation}
E(t) := \int_{\Omega}	\Big(-i\sum_{j=1}^d \Phi^*\sigma_j \partial_j \Phi + \Phi^*\sigma_3 \Phi + \varepsilon^2G(\Phi)\Big) d\bx \equiv E(0), \quad t \geq 0,
\end{equation}
where
\begin{equation}
G(\Phi) = \frac{\lambda_1}{2}\left(\Phi^*\sigma_3\Phi\right)^2	+ \frac{\lambda_2}{2}|\Phi|^4, \quad \Phi \in \mathbb{C}^2.
\end{equation}

For the NLDE \eqref{eq:NLDE_21} with $\eps = 1$, i.e., the classical regime, there are extensive analytical and numerical studies in the literature. For the existence and multiplicity of bound states and/or standing wave solutions, we refer to \citet{BCDM}, \citet{BDing}, \citet{CV}, \citet{ES}, \citet{FZ}, \citet{KA}, \citet{Pecher} and references therein. For the special case $d=1$ and $\varepsilon =1$ in the NLDE \eqref{eq:NLDE_21} with $\lambda_1 = -1$ and $\lambda_2 = 0$ in the nonlinearity \eqref{eq:nonlinear}, it admits explicit soliton solutions \citep{FS,MP,Rafflski,Taka}. In the numerical aspects, various numerical schemes have been proposed and analyzed including finite difference time domain (FDTD) methods \citep{BCJY,FY1}, exponential wave integrator Fourier pseudospectral (EWI-FP) method \citep{BCJY} and time-splitting Fourier pseudospectral (TSFP) method \citep{FSS,FLB}. However, for the NLDE \eqref{eq:NLDE_21} with $0 < \varepsilon \ll 1$, it is interesting to study the long-time dynamics for $t \in [0, T_{\eps}]$ with $T_{\varepsilon} = O(1/\eps^2)$. To the best of our knowledge, there are very few numerical analysis results on the error bounds of numerical methods for the long-time dynamics of the NLDE  \eqref{eq:NLDE_21} in the literature. Recently, we rigorously carried out the error bounds for the long-time dynamics of the Dirac equation with small potentials \citep{BFY,FY1,FY2}. Based on the results for the linear case, the TSFP method performs much better than other numerical methods with the improved uniform error bounds established in the long-time regime by employing the {\bf regularity compensation oscillation} (RCO) technique \citep{BFY}. The aim of this paper is to establish improved uniform bounds on time-splitting methods for the NLDE \eqref{eq:NLDE_21} up to the time at $O(1/\eps^2)$. Then, we combine the Fourier pseudospectral method in space and extend the improved uniform error bounds to the full-discretization. With the help of the RCO technique, the improved uniform error bounds are at $O(\eps^2\tau^2 + \tau^{m-1}_0)$ and $O(h^{m-1}+\eps^2\tau^2 + \tau^{m-1}_0)$, respectively, for the second-order semi-discretization and full-discretization for the NLDE with $O(\eps^2)$-nonlinearity up to the time at $O(1/\eps^2)$ with $m \geq 3$ depending on the regularity of the exact solution and $\tau_0 \in (0, 1)$ a fixed chosen parameter. Specifically, when the exact solution is smooth, then 
the improved uniform error bounds are at $O(\eps^2\tau^2)$ and $O(h^{m-1}+\eps^2\tau^2)$ for the second-order semi-discretization and full-discretization, respectively. They are indeed better (or improved) than
the standard error bounds at $O(\tau^2)$ and $O(h^{m-1}+\tau^2)$ in the literature for the second-order time-splitting semi-discretization and full-discretization, respectively, especially when $0<\eps\ll1$.

The main idea of the RCO technique is to choose the frequency cut-off parameter $\tau_0$ and control high frequency modes ($> 1/\tau_0$) by the regularity of the exact solution and analyze low frequency modes by phase cancellation and energy method.  Similar to  the (nonlinear) Schr\"odinger equation on an irrational rectangle and the nonlinear Klein-Gordon equation as non-periodic examples in \citep{BCF1,BCF2}, the free Dirac operator is also non-periodic in 1D. The RCO technique is not only applicable to the periodic evolutionary PDE but also to the non-periodic case.

The rest of this paper is organized as follows. In Section 2, we discretize the NLDE \eqref{eq:NLDE_21} in time by the Strang splitting method to obtain the semi-discretization and combine the Fourier pseudospectral method in space for the full-discretization. In Section 3, we establish improved uniform error bounds up to the time at $O(1/\eps^2)$ for the semi-discretization and full-discretization, respectively. In Section 4, we present some numerical results to confirm our error estimates and discuss the improved error bounds for an oscillatory NLDE. Finally, some conclusions are drawn in Section 5. Throughout this paper, we adopt the notation $A \lesssim B$ to represent that there exists a generic constant $C > 0$, which is independent of the mesh size $h$ and time step $\tau$ as well as the parameter $\varepsilon$ such that $|A| \leq C B$.
\section{Discretizations}
In this section, we apply a second-order time-splitting method to discretize the NLDE \eqref{eq:NLDE_21} and combine the Fourier pseudospectral method in space to derive the time-splitting Fourier pseudospectral (TSFP) method. For simplicity of notations, we only present the numerical methods and their analysis in 1D, i.e., $d = 1$. Numerical schemes and corresponding results can be easily generalized to the NLDE \eqref{eq:NLDE_21} in 2D and to the four-component nonlinear Dirac equation in 3D \citep{BCJY}. In 1D, the NLDE \eqref{eq:NLDE_21} on the computational domain $\Omega = (a, b)$ with periodic boundary conditions collapses to
\begin{align}
\label{eq:Dirac_1D}
&i\partial_t\Phi =  \left(- i \sigma_1 \partial_x + \sigma_3 \right)\Phi+ \varepsilon^2  {\bF}(\Phi)\Phi, \quad x \in \Omega, \quad t > 0,\\
\label{eq:ib}
&\Phi(t, a) = \Phi(t, b),\quad t \geq 0; \quad \Phi(0, x) = \Phi_0(x),\quad x \in \overline{\Omega},
\end{align}
where $\Phi := \Phi(t, x)$, $\Phi_0(a) = \Phi_0(b)$ and the nonlinearity ${\mathbf F}(\Phi)$ is given in \eqref{eq:nonlinear}.

For an integer $m \geq 0$, we denote by $H^m(\Omega)$  the set of functions $u(x) \in L^2(\Omega)$ with finite $H^m$-norm given by
\begin{equation}
\|u\|_{H^m}^2=\sum\limits_{l\in\mathbb{Z}} (1+\mu_l^2)^m\left|\widehat{u}_l\right|^2,\quad \mbox{for}\quad u(x) = \sum\limits_{l\in\mathbb{Z}}\widehat{u}_l e^{i\mu_l(x-a)},\quad \mu_l = \frac{2\pi l}{b-a},
\end{equation}
where $\widehat{u}_l (l \in \mathbb{Z})$ are the Fourier coefficients of the function $u(x)$ \citep{ST}. In fact, the space $H^m(\Omega)$ is the subspace of the classical Sobolev space $W^{m, 2}(\Omega)$, which consists of functions with derivatives of order up to $m-1$ being $(b-a)$-periodic.

Denote the index set $\mathcal{T}_M = \{l ~|~l = -M/2,-M/2+1, \ldots, M/2-1\}$ and the spaces
\begin{align*}
&X_M = \{ U = (U_0, \ldots, U_M)^T~|~U_j \in \mathbb{C}^2,\ j = 0, 1, \ldots, M,\ U_0 = U_M\},\\
& Y_M = Z_M \times Z_M, \quad Z_M = \mbox{span}\{\phi_l(x) = e^{i\mu_l(x-a)},\  l \in \mathcal{T}_M\}.
\end{align*}
The projection operator  $P_M: (L^2(\Omega))^2 \to Y_M$ is defined as
\begin{equation}
	(P_M U)(x) := \sum_{l \in \mathcal{T}_M}\widehat{U}_l e^{i\mu_l(x-a)},\quad U(x) \in (L^2(\Omega))^2,
\end{equation}
where
\begin{equation}
\widehat{U}_l = \frac{1}{b-a}\int^b_a U(x) e^{-i\mu_l(x-a)} dx, \quad l \in \mathcal{T}_M.
\end{equation}
Define the space $(C_{\rm per}(\overline{\Omega}))^2 = \{U \in (C(\overline{\Omega}))^2~|~U(a) = U(b)\}$ and the interpolation operator $I_M : (C_{\rm per}(\overline{\Omega}))^2 \to Y_M$ or $I_M : X_M \to Y_M$ as
\begin{equation}
(I_M U)(x) := \sum_{l \in \mathcal{T}_M}\widetilde{U}_l e^{i\mu_l(x-a)},\quad  U(x) \in  (C_{\rm per}(\overline{\Omega}))^2 \quad  \mbox{or} \quad  U \in X_M,
\end{equation}
where
\begin{equation}
\widetilde{U}_l = \frac{1}{M}\sum_{j=0}^{M-1}U_j e^{-2ijl\pi/M},\quad l \in \mathcal{T}_M, 	
\end{equation}
where $U_j = U(x_j)$ for a function $U(x)$.

\subsection{Semi-discretization by a second-order time-splitting method}
Denote the free Dirac operator as
\begin{equation}
\bT := -i\sigma_1\partial_x + \sigma_3,
\end{equation}
then the NLDE \eqref{eq:Dirac_1D} can be expressed as
\begin{equation}
i\partial_t\Phi(t, x) = \bT \Phi(t, x) + \eps^2  {\bf F}(\Phi(t, x))\Phi(t, x), \quad x \in \Omega, \quad t > 0.
\end{equation}

Choose the time step size as $\tau = \Delta t > 0$ and time steps $t_n = n\tau$ for $n=0,1,\ldots$. Denote by $\Phi^{[n]} := \Phi^{[n]}(x)$ the approximation of $\Phi(t_n, x)$ for $n \geq 0$, then a semi-discretization for the NLDE \eqref{eq:Dirac_1D} via the second-order time-splitting (Strang splitting)  could be expressed as \citep{BCJY,BCY2}
\begin{equation}
\Phi^{[n+1]} = \mathcal{S}_{\tau}(\Phi^{[n]}) := e^{-\frac{i\tau}{2}\bT}e^{-i\eps^2\tau{\bF}\left( e^{-\frac{i\tau}{2}\bT}\Phi^{[n]}\right)} e^{-\frac{i\tau}{2}\bT}\Phi^{[n]},\quad n = 0, 1, \ldots,
\label{eq:semi}	
\end{equation}
with the initial data taken as $\Phi^{[0]} = \Phi_0(x)$ for $x \in \overline{\Omega}$.

\begin{remark}
The second-order time-splitting (Strang splitting) method is applied to discretize the NLDE \eqref{eq:Dirac_1D} in time and it is straightforward to design the first-order Lie-Trotter splitting \citep{Trotter} and higher order schemes, e.g., the fourth-order partitioned Runge-Kutta (PRK4) splitting method \citep{BY,MQ}.
\end{remark}

\subsection{Full-discretization}
Given a spatial mesh size $h = (b - a)/M$ with $M$ an even positive integer,  the spatial grid points are $x_j := a+j h$ for $j = 0, 1, \ldots, M$. Let $\Phi^n_j$ be the numerical approximation of $\Phi(t_n, x_j)$ and denote $\Phi^n = (\Phi^n_0, \Phi^n_1, \ldots, \Phi^n_M)^T \in X_M$ as the solution vector at $t = t_n$. The initial data is taken as $\Phi^0_j = \Phi_0(x_j)$ for $j = 0, 1, \ldots, M$, then from time $t = t_n$ to $t = t_{n+1}$, the time-splitting Fourier pseudospectral (TSFP) method for discretizating the  NLDE \eqref{eq:Dirac_1D} is given as
\begin{equation}
\begin{split}
&\Phi^{(1)}_j = \sum_{l \in \mathcal{T}_M} e^{-i\frac{\tau \Gamma_l}{2}}\widetilde{(\Phi^n)}_l e^{i\mu_l(x_j-a)} =  \sum_{l \in \mathcal{T}_M} Q_l e^{-i\frac{\tau D_l}{2}}(Q_l)^T \widetilde{(\Phi^n)}_l e^{\frac{2ijl\pi}{M}},\\
& \Phi^{(2)}_j  = e^{-i\varepsilon^2 \tau {\bF}(\Phi^{(1)}_j)} \Phi^{(1)}_j = e^{-i\varepsilon^2 \tau \Lambda_j} \Phi^{(1)}_j, \quad 0 \leq j \leq M, \quad n \geq 0, \\
&\Phi^{n+1}_j = \sum_{l \in \mathcal{T}_M} e^{-i\frac{\tau \Gamma_l}{2}}\widetilde{(\Phi^{(2)})}_l e^{i\mu_l(x_j-a)} =  \sum_{l \in \mathcal{T}_M} Q_l e^{-i\frac{\tau D_l}{2}}(Q_l)^T \widetilde{(\Phi^{(2)})}_l e^{\frac{2ijl\pi}{M}},
\end{split}
\label{eq:TSFP}
\end{equation}
where $\Gamma_l = \mu_l \sigma_1 + \sigma_3 = Q_l D_l (Q_l)^T$ with $\delta_l =\sqrt{1+\mu^2_l}$,
\begin{equation}\label{eq:Qdef}
\Gamma_l = \begin{pmatrix} 1  &  \ \mu_l \\ \mu_l & \ -1 \end{pmatrix}, \quad Q_l = \begin{pmatrix} \frac{1+\delta_l}{\sqrt{2\delta_l(1+\delta_l)}}  & \  -\frac{\mu_l}{\sqrt{2\delta_l(1+\delta_l)}} \\ \frac{\mu_l}{\sqrt{2\delta_l(1+\delta_l)}} &\    \frac{1+\delta_l}{\sqrt{2\delta_l(1+\delta_l)}} \end{pmatrix}, \quad D_l = \begin{pmatrix} \delta_l  & \  0 \\ 0 &\  -\delta_l \end{pmatrix},
\end{equation}
and $\Lambda_j = \textrm{diag}(\Lambda_{j, +},\Lambda_{j, -})$ with $\Lambda_{j, \pm} = \lambda_2 |\Phi^{(1)}_j|^2 \pm  \lambda_1(\Phi^{(1)}_j)^{\ast}\sigma_3\Phi^{(1)}_j$.

\section{Improved uniform error bounds}
In this section, we rigorously prove the improved uniform error bounds for the second-order time-splitting method in propagating the NLDE with $O(\eps^2)$ nonlinearity up to the long-time at $O(1/\eps^2)$. For the simplicity of presentation, we shall assume $\lambda_1=0$ in the subsequent discussion, where the results and the proof are also valid if $\lambda_1\neq0$ by using the same arguments.

\subsection{Main results}
We assume the exact solution $\Phi(t, x)$ of the NLDE \eqref{eq:Dirac_1D} up to the time at $T_{\eps} = T/\eps^2$ with $T  > 0$ fixed satisfies
\begin{equation*}
\textrm{(A)} \hspace{4cm}  \Phi \in  {L^{\infty}([0, T_\eps]; (H^{m}(\Omega))^2)}, \quad m \geq 3. \hspace{4cm}
\end{equation*}
Let $\Phi^{[n]}$ be the approximation obtained from the time-splitting method \eqref{eq:semi} and $\lambda_1 = 0$ in \eqref{eq:nonlinear}. According to the standard analysis in \citet{BCJY}, under the assumption (A), for sufficiently small $0 < \tau \leq \tau_c$ with $\tau_c > 0$ a constant, there exists a constant $M > 0$ depending on $T$ and $\left\|\Phi\right\|_{L^{\infty}([0, T_\eps]; (H^{m})^2)}$ such that
\begin{equation}
\left\|\Phi^{[n]}\right\|_{H^m} \leq M, \quad 0 \leq n \leq \frac{T/\eps^2}{\tau}.	
\label{eq:nb}
\end{equation}
In this work, we will establish the following improved uniform error bounds up to the long-time $T_\eps$.

\begin{theorem}
Under the assumption (A), for $0 < \tau_0 < 1$ sufficiently small and independent of $\eps$ such that, when $0 < \tau \leq \alpha \frac{\pi (b-a)\tau_0}{2\sqrt{\tau^2_0(b-a)^2+4\pi^2(1+\tau^2_0)}}$ for a fixed constant $\alpha \in (0, 1)$, we have the following improved uniform error bound for any $\eps \in (0, 1]$
\begin{equation}
\left\|\Phi(t_n, x) - \Phi^{[n]}\right\|_{H^1} \lesssim \eps^2\tau^2 + \tau_0^{m-1},\quad 0 \leq n \leq \frac{T/\eps^2}{\tau}.
\label{eq:semi_bound}	
\end{equation}
In particular, if the exact solution is sufficiently smooth, e.g. $\Phi(t, x) \in (H^{\infty}(\Omega))^2$, the last term $\tau_0^{m-1}$ decays exponentially fast and could be ignored practically for small enough $\tau_0$, and the improve uniform error bound would become
\begin{equation}
\left\|\Phi(t_n, x) - \Phi^{[n]}\right\|_{H^1} \lesssim \eps^2\tau^2, \quad 0 \leq n \leq \frac{T/\eps^2}{\tau}.
\end{equation}
\label{thm:semi}
\end{theorem}

Let $\Phi^n$ be the numerical approximation obtained from the TSFP \eqref{eq:TSFP}, then we have the following improved uniform error bounds for the full-discretization.

\begin{theorem}
Under the assumption (A), there exist $h_0 >0$ and $0 < \tau_0 < 1$ sufficiently small and independent of $\eps$ such that for any $0 < \eps \leq 1$, when $0 < h \leq h_0$ and $0 < \tau \leq \alpha \frac{\pi (b-a)\tau_0}{2\sqrt{\tau^2_0(b-a)^2+4\pi^2(1+\tau^2_0)}}$ for a fixed constant $\alpha \in (0, 1)$, the following improved uniform error bound holds
\begin{equation}
\left\|\Phi(t_n, x) - I_M\Phi^{n}\right\|_{H^1} \lesssim h^{m-1} + \eps^2\tau^2 + \tau_0^{m-1},\quad 0 \leq n \leq \frac{T/\eps^2}{\tau}.
\label{eq:full_bound}	
\end{equation}
Similarly, for the sufficiently smooth exact solution, e.g. $\Phi(t, x) \in (H^{\infty}(\Omega))^2$, the improve uniform error bound would become
\begin{equation}
\left\|\Phi(t_n, x) - I_M\Phi^{n}\right\|_{H^1} \lesssim h^{m-1} + \eps^2\tau^2, \quad 0 \leq n \leq \frac{T/\eps^2}{\tau}.
\end{equation}
\label{thm:TSFP}
\end{theorem}

\begin{remark}
Here, we prove $H^1$-error bounds for 1D problem to control the nonlinearity since $H^1(\mathbb{R})$ is an algebra.  Corresponding error estimates should be in $H^2$-norm for 2D and 3D cases (2D case in the sense of \eqref{eq:NLDE_21}, and 3D case in the sense of the four-component NLDE given in \citep{BCJY}). Of course, higher regularity assumptions of the exact solution are required for higher order Sobolev norm estimates.
\end{remark}

\begin{remark} Under appropriate assumptions of the exact solution, the improved uniform error bounds could be extended to the first-order Lie-Trotter splitting and the fourth-order PRK splitting method with improved uniform error bounds at $\varepsilon^2 \tau$ and $\varepsilon^2 \tau^{4}$, respectively.
\end{remark}

\begin{remark}
For the NLDE \eqref{eq:Dirac_1D} with general matrix nonlinearity $\bF(\Phi)$ in \eqref{eq:nonlinear} when $\lambda_1 \neq 0$, the proof is similar and we omit the details here for brevity.	
\end{remark}

\begin{remark}
Define the discrete energy at $t = t_n$ with the mesh size $h$ as
\begin{equation}
E^n_h = h\sum_{j=0}^{M-1}\left[-i(\Phi_j^n)^*\sigma_1 (\Phi')^n_j + (\Phi_j^n)^*\sigma_3 \Phi^n_j + \varepsilon^2G(\Phi^n_j)\right],
\end{equation}
where
\begin{equation}
(\Phi')^n_j = i\sum_{l \in \mathcal{T}_M}\mu_l \widetilde{(\Phi^n)}_l e^{i\mu_l(x_j-a)}, \quad j = 0, 1, \ldots, M-1,
\label{eq:prime}
\end{equation}
then we have the following estimate for the discrete energy
\begin{equation}
\left|E^n_h - E^0_h\right|\lesssim h^{m-1}+\eps^2\tau^2 + \tau_0^{m-1}, \quad 0 \leq n \leq \frac{T/\eps^2}{\tau}.
\end{equation}
Similarly, for the sufficiently smooth exact solution, e.g. $\Phi(t, x) \in (H^{\infty}(\Omega))^2$, the estimate for the discrete energy would become
\begin{equation}
\left|E^n_h - E^0_h\right|\lesssim h^{m-1}+\eps^2\tau^2, \quad 0 \leq n \leq \frac{T/\eps^2}{\tau}.
\end{equation}
\end{remark}

\subsection{Proof for Theorem \ref{thm:semi}}
Denote the exact solution flow $\Phi(t_n) \to \Phi(t_{n+1})$ as
\begin{equation}
\Phi(t_{n+1}) = \mathcal{S}_{e, \tau}(\Phi(t_n)),\quad 0 \leq n \leq \frac{T/\varepsilon^2}{\tau},
\end{equation}
where we take $\Phi(t_n) := \Phi(t_n, x)$ for simplicity. By the standard analysis, we have the following estimates for the local truncation error \citep{BCJY}.
\begin{lemma}
For $0 < \eps \leq 1$, the local truncation error for the discrete-in-time second-order splitting \eqref{eq:semi} could be written as
\begin{equation}
\mathcal{E}^n := \mathcal{S}_{\tau}(\Phi(t_{n})) - \Phi(t_{n+1}) = \mathcal{F}(\Phi(t_n)) + \mathcal{R}^n, \quad n = 0, 1, \ldots,
\label{eq:local-dec}
\end{equation}	
where
\begin{align}
&\mathcal{F}(\Phi(t_n)) = i\eps^2 \left(\int^{\tau}_0 f_s(\Phi(t_n)) ds - \tau f_{\tau/2}(\Phi(t_n))\right),	 \label{eq:matF}\\
&f_s(\Phi(t_n))=e^{-i(\tau-s) \bT}\left({\mathbf F}(e^{-is\bT}\Phi(t_n)) e^{-is\bT}\Phi(t_{n})\right)\label{eq:fndef},
\end{align}
and under the assumption (A), the following error bounds hold
\begin{equation}
\left\|\mathcal{F}(\Phi(t_n))\right\|_{H^1} \lesssim \eps^2\tau^3,	\quad \left\|\mathcal{R}^n\right\|_{H^1} \lesssim \eps^4\tau^3.
\label{eq:lb}
\end{equation}
\end{lemma}

In addition, under the assumption (A), the nonlinearity could be controlled by the estimates \eqref{eq:nb} for the numerical solution $\Phi^{[n]}$. Introduce the error function
\begin{equation}
{\bf e}^{[n]} := {\bf e}^{[n]}(x)  = \Phi^{[n]} - \Phi(t_n),	\quad n= 0, 1, \ldots,
\end{equation}
 we have the error equation from \eqref{eq:semi} and \eqref{eq:local-dec} as
\begin{equation}
{\bf e}^{[n+1]} = \mathcal{S}_{\tau}(\Phi^{[n]})	 - \mathcal{S}_{\tau}(\Phi(t_n)) + \mathcal{E}^n = e^{-i\tau\bT}{\bf e}^{[n]} + W^n + \mathcal{E}^n, \quad n \geq 0,
\label{eq:egrow}
\end{equation}
where $W^n := W^n(x)$ is given by
\begin{equation*}
W^n(x) = e^{-\frac{i\tau}{2}\bT}\left[\left(e^{-i\eps^2\tau\bF(e^{-\frac{i\tau}{2}\bT} \Phi^{[n]})}-I_2\right)e^{-\frac{i\tau}{2}\bT} \Phi^{[n]} - \left(e^{-i\eps^2\tau\bF(e^{-\frac{i\tau}{2}\bT} \Phi(t_n))}-I_2\right)e^{-\frac{i\tau}{2}\bT} \Phi(t_n)\right],
\end{equation*}
with the bound implied by \eqref{eq:nb}
\begin{equation}
\left\|W^n(x)\right\|_{H^1} \lesssim \eps^2 \tau\left\|{\bf e}^{[n]}\right\|_{H^1}.	
\label{eq:wb}
\end{equation}
Based on \eqref{eq:egrow}, we obtain
\begin{equation}
{\bf e}^{[n+1]} = 	e^{-i(n+1)\tau\bT}{\bf e}^{[0]} + \sum\limits_{k=0}^n 	e^{-i(n-k)\tau\bT}\left(W^k + \mathcal{E}^k\right), \quad 0 \leq n \leq \frac{T/\eps^2}{\tau} - 1.
\end{equation}
Noticing ${\bf e}^{[0]} = {\bf 0}$, combining \eqref{eq:local-dec}, \eqref{eq:lb} and \eqref{eq:wb}, we have the estimates for $0 \leq n \leq T_\eps/\tau-1$,
\begin{equation}
\left\|{\bf e}^{[n+1]}\right\|_{H^1} \lesssim \eps^2\tau^2 + \eps^2\tau \sum\limits_{k=0}^n	\left\|{\bf e}^{[k]}\right\|_{H^1} + \left\|\sum\limits_{k=0}^n e^{-i(n-k)\tau\bT}\mathcal{F}(\Phi(t_k))\right\|_{H^1}.
\label{eq:semi_finial}
\end{equation}
In order to obtain the improved uniform error bounds \eqref{eq:semi_bound},  we will employ the {\textbf{regularity compensation oscillation}} (RCO) technique \citep{BCF1,BCF2} to deal with the last term on the RHS of \eqref{eq:semi_finial}. From the NLDE \eqref{eq:Dirac_1D}, we find that $\partial_t\Phi + i\bT\Phi = O(\eps^2)$.  It is natural to consider the `twisted variable'
\begin{equation}
\Psi(t, x) = e^{it\bT}\Phi(t, x), \quad t \geq 0,	
\end{equation}
which satisfies the equation $\partial_t \Psi(t, x) = \eps^2 e^{it\bT}\left(\bF\left(e^{-it\bT}\Psi(t, x)\right)e^{-it\bT}\Psi(t, x)\right)$.  Under the assumption (A), we have $\|\Psi\|_{L^{\infty}([0, T_\eps]; (H^m(\Omega))^2)}	\lesssim 1$ and $
\|\partial_{t}\Psi\|_{L^{\infty}([0, T_\eps]; (H^m(\Omega))^2)} \lesssim \eps^2$ with
\begin{equation}
\left\|\Psi(t_{n+1}) - \Psi(t_n)\right\|_{H^m}	\lesssim \eps^2 \tau, \quad 0 \leq n \leq \frac{T/\eps^2}{\tau} - 1.
\label{eq:S1twist}
\end{equation}

{\textbf{Step 1}}. Choose the cut-off parameter on the Fourier modes. Let $\tau_0 \in (0, 1)$
  and $M_0 = 2\lceil 1/\tau_0 \rceil \in \mathbb{Z}^+$ ($\lceil \cdot \rceil$ is the ceiling function) with $1/\tau_0 \leq M_0/2 < 1+ 1/\tau_0$. Under the assumption (A), we have the following estimate
\begin{equation}
\left\|P_{M_0}\mathcal{F}(e^{-it_k\bT}P_{M_0}\Psi(t_k)) - \mathcal{F}(e^{-it_k\bT}\Psi(t_k))\right\|_{H^1}	\lesssim \eps^2 \tau \tau_0^{m-1}.
\end{equation}
Based on above estimates, \eqref{eq:semi_finial} would imply for $0 \leq n \leq T_\eps/\tau-1$,
\begin{equation}
\left\|{\bf e}^{[n+1]}\right\|_{H^1}\lesssim \tau_0^{m-1} + \eps^2\tau^2 + \eps^2\tau \sum_{k=0}^{n} \left\|{\bf e}^{[k]}\right\|_{H^1} + \left\| \mathcal{L}^n\right\|_{H^1},
\label{eq:finial1}
\end{equation}
with
\begin{equation}
\mathcal{L}^n = \sum_{k=0}^{n} e^{i(k+1)\tau \bT}P_{M_0}\mathcal{F}\left(e^{-ik\tau\bT}(P_{M_0}\Psi(t_k))\right).
\label{eq:matl}
\end{equation}

{\textbf{Step 2}}. Analyze the low Fourier modes term $\mathcal{L}^n$. For $l\in\mathcal{T}_{M_0}$, define the index set $\mathcal{I}_l^{M_0}$ associated to $l$ as
\begin{equation}
\mathcal{I}_l^{M_0}=\left\{(l_1,l_2,l_3)\ | \ l_1- l_2 + l_3 = l,\ l_1, l_2, l_3\in\mathcal{T}_{M_0}\right\}.
\end{equation}
We introduce
\begin{equation}
\Pi_l^+=Q_l\text{diag}(1,0)(Q_l)^T,\quad \Pi_l^-=Q_l\text{diag}(0,1)(Q_l)^T,
\end{equation}
where $\Pi_l^{\pm}$ are the projectors onto the eigenspaces of $\Gamma_l$ corresponding to the eigenvalues $\pm\delta_l$, respectively. Moreover, we have $(\Pi_l^\pm)^T=\Pi_l^\pm$, $\Pi_l^++\Pi_l^-=I_2$, $(\Pi_l^\pm)^2=\Pi_l^\pm$, $\Pi_l^{\pm}\Pi_l^{\mp}={\bf 0}$.
By direct computation, we  have
\begin{equation}
e^{it{\bT}}P_{M_0}\Psi(t_k)=\sum\limits_{l\in\mathcal{T}_{M_0}}\left(e^{it\delta_l}\Pi_l^++e^{-it\delta_l}\Pi_l^-\right)\widehat{\Psi}_l(t_k)e^{i\mu_l(x-a)}.
\end{equation}
In view of the definition of $\mathcal{F}$ \eqref{eq:matF} , we have the following expansion
\begin{align*}
e^{i(k+1)\tau\bT}P_{M_0}(f_s(e^{-i k\tau\bT}P_{M_0}\Psi(t_k))) 
=\lambda_2\sum\limits_{l\in\mathcal{T}_{M_0}}\sum\limits_{(l_1,l_2,l_3)\in\mathcal{I}_l^{M_0}} \sum\limits_{\nu_j=\pm,j=1,2,3,4}\mathcal{G}^{\nu_1, \nu_2, \nu_3, \nu_4}_{k, l, l_1, l_2, l_3}(s) e^{i\mu_l(x-a)},
\end{align*}
where the coefficients  $\mathcal{G}^{\nu_1, \nu_2, \nu_3, \nu_4}_{k, l, l_1, l_2, l_3}(s)$ are functions of $s$ defined as
\begin{equation}
\mathcal{G}^{\nu_1, \nu_2, \nu_3, \nu_4}_{k, l, l_1, l_2, l_3}(s)
= e^{i(t_k+s)\delta^{\nu_1, \nu_2, \nu_3, \nu_4}_{l, l_1, l_2, l_3}}\Pi_l^{\nu_1}(\widehat{\Psi}_{l_2}(t_k))^*\Pi_{l_2}^{\nu_3}\Pi_{l_1}^{\nu_2}\widehat{\Psi}_{l_1}(t_k)\Pi_{l_3}^{\nu_4} \widehat{\Psi}_{l_3}(t_k)
\label{eq:mGdef}
\end{equation}
with $\delta^{\nu_1, \nu_2, \nu_3, \nu_4}_{l, l_1, l_2, l_3} = \nu_1 \delta_l-\nu_2 \delta_{l_1}+\nu_3 \delta_{l_2}-\nu_4 \delta_{l_3}$. Thus, we have
\begin{equation}
\mathcal{L}^n = \lambda_2i\eps^2\sum\limits_{k=0}^n
\sum\limits_{l\in\mathcal{T}_{M_0}}\sum\limits_{(l_1, l_2, l_3)\in\mathcal{I}_l^{M_0}}\sum\limits_{\nu_j=\pm, j = 1, 2, 3, 4}\Upsilon_{k, l, l_1, l_2, l_3}^{\nu_1, \nu_2, \nu_3, \nu_4}e^{i\mu_l(x - a)},
\label{eq:reminder-def}
\end{equation}
where
\begin{align}
\Upsilon_{k, l, l_1, l_2, l_3}^{\nu_1, \nu_2, \nu_3, \nu_4}  = -\tau \mathcal{G}^{\nu_1, \nu_2, \nu_3, \nu_4}_{k, l, l_1, l_2, l_3}(\tau/2) + \int^{\tau}_0 \mathcal{G}^{\nu_1, \nu_2, \nu_3, \nu_4}_{k, l, l_1, l_2, l_3}(s) \, ds	
 =  r_{l, l_1, l_2, l_3}^{\nu_1, \nu_2, \nu_3, \nu_4}e^{it_k\delta^{\nu_1, \nu_2, \nu_3, \nu_4}_{l, l_1, l_2, l_3}}c_{k, l, l_1, l_2, l_3}^{\nu_1, \nu_2, \nu_3, \nu_4},
\label{eq:Lam}
\end{align}
and the vector coefficients $c_{l, l_1, l_2, l_3}^{\nu_1, \nu_2, \nu_3, \nu_4}$ and the scalar coefficient $r_{l, l_1, l_2, l_3}^{\nu_1, \nu_2, \nu_3, \nu_4}$ are given as
\begin{align}
c_{l, l_1, l_2, l_3}^{\nu_1, \nu_2, \nu_3, \nu_4}& =
\Pi_l^{\nu_1}(\widehat{\Psi}_{l_2}(t_k))^*\Pi_{l_2}^{\nu_3}
\Pi_{l_1}^{\nu_2}\widehat{\Psi}_{l_1}(t_k)\Pi_{l_3}^{\nu_4} \widehat{\Psi}_{l_3}(t_k),\label{eq:calFdef}\\
r_{l, l_1, l_2, l_3}^{\nu_1, \nu_2, \nu_3, \nu_4}& = -\tau e^{i\tau\delta_{k, l, l_1, l_2, l_3}^{\nu_1, \nu_2, \nu_3, \nu_4}/2}+\int_0^\tau e^{is\delta_{k, l, l_1, l_2, l_3}^{\nu_1, \nu_2, \nu_3, \nu_4}}\,ds \nn \\
& =O\left(\tau^3 (\delta_{k, l, l_1, l_2, l_3}^{\nu_1, \nu_2, \nu_3, \nu_4})^2\right).\label{eq:rest}
\end{align}
We only need to consider the case $\delta_{k, l, l_1, l_2, l_3}^{\nu_1, \nu_2, \nu_3, \nu_4} \neq 0$ as $r_{l, l_1, l_2, l_3}^{\nu_1, \nu_2, \nu_3, \nu_4} = 0$ if  $\delta_{k, l, l_1, l_2, l_3}^{\nu_1, \nu_2, \nu_3, \nu_4}= 0$. For $l \in \mathcal{T}_{M_0}$ and $(l_1, l_2, l_3) \in \mathcal{I}^{M_0}_l$, we have
\begin{equation}
|\delta_{l, l_1, l_2, l_3}^{\nu_1, \nu_2, \nu_3, \nu_4}| \leq 4\delta_{M_0/2}=4\sqrt{1+\mu_{M_0/2}^2} < 4\sqrt{1+\frac{4\pi^2(1+\tau_0)^2}{\tau_0^2(b - a)^2}},
\end{equation}
which implies  $\frac{\tau}{2} |\delta_{l, l_1, l_2, l_3}^{\nu_1, \nu_2, \nu_3, \nu_4}| < \alpha \pi$ when $0 < \tau \leq \alpha \frac{\pi(b-a)\tau_0}{2\sqrt{\tau_0^2(b-a)^2+4\pi^2(1+\tau_0)^2}} := \tau^{\alpha}_0$ $ (0<\tau_0, \alpha < 1)$. Since $\mathcal{G}^{\nu_1, \nu_2, \nu_3, \nu_4}_{k, l, l_1, l_2, l_3}$ \eqref{eq:mGdef} is similar, it suffice to consider the typical case $\nu_1 = \nu_2 = \nu_3 = \nu_4 = +$ . Denoting $S^n_{l, l_1, l_2, l_3} = \sum_{k=0}^ne^{it_k\delta_{l, l_1, l_2, l_3}^{+,+,+,+}}$ ($ n \geq 0$) ,  for $0 < \tau < \tau^{\alpha}_0$,
we have
\begin{equation}
\label{eq:Sbd}
|S^n_{l, l_1, l_2, l_3}| \leq  \frac{1}{|\sin(\tau \delta_{l,l_1,l_2,l_3}^{+,+,+,+}/2)|}\leq\frac{C}{\tau|\delta_{l,l_1,l_2,l_3}^{+,+,+,+}|},\quad\forall n\ge0,
\end{equation}
with $C = \frac{2\alpha\pi}{\sin(\alpha\pi)}$. Using summation-by-parts, we find from \eqref{eq:Lam} that
\begin{equation}
\sum_{k=0}^n \Upsilon_{k,l,l_1,l_2,l_3}^{+,+,+,+}=
r_{l,l_1,l_2,l_3,l_4}^{+,+,+,+}\left[\sum_{k=0}^{n-1}S^k_{l, l_1, l_2, l_3} \left(c_{k,l,l_1,l_2,l_3}^{+,+,+,+}-c_{k+1,l,l_1,l_2,l_3}^{+,+,+,+}\right) +S^n_{l, l_1, l_2, l_3} c_{n,l,l_1,l_2,l_3}^{+,+,+,+}\right],
\label{eq:lambdasum}
\end{equation}
with
\begin{align}
c_{k,l,l_1,l_2,l_3}^{+,+,+,+} - c_{k+1,l,l_1,l_2,l_3}^{+,+,+,+} = & \ \Pi_l^{+}\left((\widehat{\Psi}_{l_2}(t_k))^*\Pi_{l_2}^{+}\Pi_{l_1}^{+}(\widehat{\Psi}_{l_1}(t_k)-\widehat{\Psi}_{l_1}(t_{k+1}))\right)\Pi_{l_3}^{+} \widehat{\Psi}_{l_3}(t_k) \nonumber\\
&+\Pi_l^{+}\left(\widehat{\Psi}_{l_2}(t_k)-\widehat{\Psi}_{l_2}(t_{k+1})\right)^*\Pi_{l_2}^{+}\Pi_{l_1}^{+}\widehat{\Psi}_{l_1}(t_{k+1})\Pi_{l_3}^{+} \widehat{\Psi}_{l_3}(t_k)\nonumber\\
&+\Pi_l^{+}(\widehat{\Psi}_{l_2}(t_{k+1}))^*\Pi_{l_2}^{+}\Pi_{l_1}^{+}\widehat{\Psi}_{l_1}(t_{k+1})\Pi_{l_3}^{+} \left(\widehat{\Psi}_{l_3}(t_k)-\widehat{\Psi}_{l_3}(t_{k+1})\right).\label{eq:cksum}
\end{align}
Combining \eqref{eq:rest}, \eqref{eq:Sbd}, \eqref{eq:lambdasum} and \eqref{eq:cksum}, we have
\begin{align}
\left|\sum_{k=0}^n\Upsilon_{k,l,l_1,l_2,l_3}^{+,+,+,+}\right|
\lesssim &\ \tau^2|\delta_{l,l_1,l_2,l_3}^{+,+,+,+}|\sum\limits_{k=0}^{n-1}\bigg(
\left|\widehat{\Psi}_{l_1}(t_k)-\widehat{\Psi}_{l_1}(t_{k+1})\right|\left|\widehat{\Psi}_{l_2}(t_k)\right| \left|\widehat{\Psi}_{l_3}(t_k)\right|\nonumber\\
&+\left|\widehat{\Psi}_{l_1}(t_{k+1})\right|\left|\widehat{\Psi}_{l_2}(t_k)-\widehat{\Psi}_{l_2}(t_{k+1})\right| \left|\widehat{\Psi}_{l_3}(t_k)\right|\nonumber\\
&+\left|\widehat{\Psi}_{l_1}(t_{k+1})\right|\left|\widehat{\Psi}_{l_2}(t_{k+1})\right| \left|\widehat{\Psi}_{l_3}(t_k)-\widehat{\Psi}_{l_3}(t_{k+1})\right|\bigg)\nonumber\\
&+ \tau^2|\delta_{l,l_1,l_2,l_3}^{+,+,+,+}| \left|\widehat{\Psi}_{l_1}(t_n)\right|\left|\widehat{\Psi}_{l_2}(t_n)\right| \left|\widehat{\Psi}_{l_3}(t_n)\right|.\label{eq:sumlambda}
\end{align}
The same estimates \eqref{eq:sumlambda} above hold for $\sum_{k=0}^n\Upsilon_{k,l,l_1,l_2,l_3}^{\pm,\pm,\pm,\pm}$ ($l\in\mathcal{T}_{M_0},\,(l_1,l_2,l_3)\in\mathcal{T}_l^{M_0}$). For the $H^1$-norm estimates, we notice that for $l\in\mathcal{T}_{M_0}$ and $(l_1,l_2,l_3)\in\mathcal{I}_l^{M_0}$, there holds
\begin{equation}
\left(1+|\mu_l|\right)|\delta_{l, l_1, l_2, l_3}^{\nu_1, \nu_2, \nu_3, \nu_4}|\lesssim \prod_{j=1}^3\left(1+\mu_{l_j}^2\right).
\label{eq:mlbd}
\end{equation}
Based on \eqref{eq:reminder-def}, \eqref{eq:sumlambda} and \eqref{eq:mlbd}, we have
\begin{align}\label{eq:sumlambda-2}
\left\|\mathcal{L}^n\right\|^2_{H^1}=&\lambda_2^2\varepsilon^4
\sum\limits_{l\in\mathcal{T}_{M_0}}\left(1+\mu_l^2\right)\left|\sum\limits_{(l_1,l_2,l_3)\in\mathcal{I}_l^{M_0}}\sum\limits_{\nu_j=\pm}\sum\limits_{k=0}^n\Upsilon_{k, l, l_1, l_2, l_3}^{\nu_1, \nu_2, \nu_3, \nu_4}\right|^2 \nn \\
\lesssim&\ \eps^4\tau^4
\bigg\{\sum_{l\in\mathcal{T}_{M_0}}\bigg(\sum\limits_{(l_1,l_2,l_3)\in\mathcal{I}_l^{M_0}}\left|\widehat{\Psi}_{l_1}(t_n)\right|\left|\widehat{\Psi}_{l_2}(t_n)\right| \left|\widehat{\Psi}_{l_3}(t_n)\right|\prod_{j=1}^3(1+\mu_{l_j}^2) \bigg)^2 \nn\\
&+n\sum\limits_{k=1}^{n-1} \bigg[\sum_{l\in\mathcal{T}_{M_0}}\bigg(\sum\limits_{(l_1,l_2,l_3)\in\mathcal{I}_l^{M_0}}
\left|\widehat{\Psi}_{l_1}(t_k)-\widehat{\Psi}_{l_1}(t_{k+1})\right|\left|\widehat{\Psi}_{l_2}(t_k)\right| \left|\widehat{\Psi}_{l_3}(t_k)\right|\prod_{j=1}^3(1+\mu_{l_j}^2) \bigg)^2\nonumber\\
&+\bigg(\sum\limits_{(l_1,l_2,l_3)\in\mathcal{I}_l^{M_0}}
\left|\widehat{\Psi}_{l_1}(t_{k+1})\right|\left|\widehat{\Psi}_{l_2}(t_k)-\widehat{\Psi}_{l_2}(t_{k+1})\right| \left|\widehat{\Psi}_{l_3}(t_k)\right|\prod_{j=1}^3(1+\mu_{l_j}^2) \bigg)^2\nonumber\\
&+\bigg(\sum\limits_{(l_1,l_2,l_3)\in\mathcal{I}_l^{M_0}}
\left|\widehat{\Psi}_{l_1}(t_{k+1})\right|\left|\widehat{\Psi}_{l_2}(t_{k+1})\right| \left|\widehat{\Psi}_{l_3}(t_k)-\widehat{\Psi}_{l_3}(t_{k+1})\right|\prod_{j=1}^3(1+\mu_{l_j}^2) \bigg)^2\bigg]\bigg\}.
\end{align}
Introduce the auxiliary function $\Theta(x)=\sum_{l\in\mathbb{Z}}(1+\mu_l^2)\left|\widehat{\Psi}_l(t_n)\right|e^{i\mu_l(x-a)}$, where $\Theta(x)\in  H^{m-2}$ implied by assumption (A) and $\|\Theta\|_{H^s}\lesssim\|\Psi(t_n)\|_{H^{s+2}}$ ($s\leq m-2$). Expanding
\begin{equation*}
|\Theta(x)|^2\Theta(x)=\sum\limits_{l\in\mathcal{Z}}\sum\limits_{l_1-l_2+l_3=l, l_j\in\mathbb{Z}} \prod_{j=1}^3\left((1+\mu_{l_j}^2)\left|\widehat{\Psi}_{l_j}(t_n)\right|\right) e^{i\mu_l(x-a)},
\end{equation*}
we could obtain that
\begin{align}
&\sum_{l\in\mathcal{T}_{M_0}}\bigg(\sum\limits_{(l_1,l_2,l_3)\in\mathcal{I}_l^{M_0}}\left|\widehat{\Psi}_{l_1}(t_n)\right|\left|\widehat{\Psi}_{l_2}(t_n)\right| \left|\widehat{\Psi}_{l_3}(t_n)\right|\prod_{j=1}^3(1+\mu_{l_j}^2) \bigg)^2\nonumber\\
&\leq\sum_{l\in\mathbb{Z}}\bigg(\sum\limits_{l_1-l_2+l_3= l, l_j\in\mathbb{Z}}\left|\widehat{\Psi}_{l_1}(t_n)\right|\left|\widehat{\Psi}_{l_2}(t_n)\right| \left|\widehat{\Psi}_{l_3}(t_n)\right|\prod_{j=1}^3(1+\mu_{l_j}^2) \bigg)^2\nonumber\\
&=\||\Theta(x)|^2\Theta(x)\|_{L^2}^2\lesssim \|\Theta(x)\|_{H^1}^6\lesssim \|\Psi(t_n)\|_{H^3}^6 \lesssim1.
\end{align}
Thus, in light of \eqref{eq:S1twist}, we could obtain the estimate for each term in \eqref{eq:matl} similarly as
\begin{align}
\left\| \mathcal{L}^n\right\|_{H^1}^2 & \lesssim \eps^4\tau^4 \bigg[\|\Psi(t_k)\|_{H^3}^6+n\sum\limits_{k=1}^{n-1}
\|\Psi(t_k)-\Psi(t_{k+1})\|^2_{H^3}(\|\Psi(t_k)\|_{H^3}+\|\Psi(t_{k+1})\|_{H^3})^4\bigg]\nn \\
& \lesssim \eps^4\tau^4+n\eps^4\tau^4 (\eps^2\tau)^2
\lesssim \eps^4\tau^4,\quad 0 \leq n\leq\frac{T/\varepsilon^2}{\tau}-1.
\label{eq:est-h1}
\end{align}
Finally, combining \eqref{eq:finial1} and  \eqref{eq:est-h1}, we have
\begin{equation}
\|{\bf e}^{[n+1]}\|_{H^1}\lesssim \tau_0^{m-1}+\eps^2\tau^2+\varepsilon^2\tau\sum_{k=0}^n\|{\bf e}^{[k]}\|_{H^1},\quad 0\leq n\leq\frac{T/\varepsilon^2}{\tau}-1.
\end{equation}
Discrete Gronwall inequality would yield
\begin{equation}
\|{\bf e}^{[n+1]}\|_{H^1}\lesssim \tau_0^{m-1}+\eps^2\tau^2, \quad 0 \leq n\leq\frac{T/\varepsilon^2}{\tau}-1,
\end{equation}
and
the improved uniform error bound \eqref{eq:semi_bound} in Theorem \ref{thm:semi} holds.

\subsection{Proof for Theorem \ref{thm:TSFP}}
Following the standard analysis in \citet{BCJY}, under the assumption (A), for $0 < \tau \leq \tau_c$ and $0 < h \leq h_c$ ($\tau_c$ and $h_c$ are two constants independent of $\eps$), there exists a constant $B$ such that the numerical solution satisfies $\left\|I_M\Phi^n\right\|_{H^1} \leq B$ for $0 \leq n \leq \frac{T/\eps^2}{\tau}$. Then we are going to establish the following improved uniform error bounds for the full-discretization.

 From the error estimates in the previous section, we have for $0 \leq n\leq \frac{T/\eps^2}{\tau}$,
\begin{equation}
\left\|\Phi(t_n, x) - \Phi^{[n]}\right\|_{H^1} \lesssim \eps^2\tau^2 + \tau^{m-1}_0,\quad	\left\|\Phi^{[n]}-P_M\Phi^{[n]}\right\|_{H^1} \lesssim h^{m-1}.
\end{equation}
Since $\Phi(t_n, x) - I_M\Phi^n = \Phi(t_n, x) - \Phi^{[n]}+ \Phi^{[n]} - P_M\Phi^{[n]} + P_M\Phi^{[n]} - I_M\Phi^n$, we derive
\begin{equation}
\left\|\Phi(t_n, x) - I_M\Phi^n\right\|_{H^1} \leq \left\|P_M\Phi^{[n]} - I_M\Phi^n\right\|_{H^1} + C_1\left(h^{m-1} + \eps^2\tau^2 + \tau^{m-1}_0 \right),
\label{eq:diff}
\end{equation}
where $C_1$ is a constant independent of $h$, $\tau$, $n$, $\eps$ and $\tau_0$. As a result, it remains to establish the estimates  on the error function ${\widetilde{\bf e}}^n := {\widetilde{\bf e}}^n(x) \in Y_M$  given as
\[
{\widetilde{\bf e}}^n := P_M\Phi^{[n]} - I_M\Phi^n,\quad 0 \leq n \leq \frac{T/\eps^2}{\tau}.\]
From \eqref{eq:TSFP}, we have
\begin{align*}
& I_M\Phi^{n+1} = e^{-\frac{i\tau}{2}{\bT}}\left(I_M\left(e^{-i\eps^2\tau{\bF}(e^{-\frac{i\tau}{2}{\bT}}I_M\Phi^n)}  e^{-\frac{i\tau}{2}{\bT}}I_M\Phi^n\right)\right),\\
& P_M\Phi^{[n+1]} =  e^{-\frac{i\tau}{2} {\bT}}\left(P_M\left(e^{-i\eps^2\tau{\bF}(e^{-\frac{i\tau}{2} {\bT}}\Phi^{[n]})} e^{-\frac{i\tau}{2} {\bT}}\Phi^{[n]}\right)\right),
\end{align*}
which lead to
\begin{align}\label{eq:errg:f}
\widetilde{\bf e}^{n+1}=e^{-i\tau\bT}\widetilde{\bf e}^{n}+Q^n(x),
\end{align}
where
\begin{align*}
Q^n(x) = & \  e^{-\frac{i\tau }{2}{\bT}}I_M\left((e^{-i\eps^2\tau{\bF}(e^{-\frac{i\tau }{2}{\bT}}I_M\Phi^n)}-I_2)e^{-\frac{i\tau }{2}{\bT}}I_M\Phi^n\right) \nn \\
& \ - e^{-\frac{i\tau }{2}{\bT}} P_M\left((e^{-i\eps^2\tau{\bF}(e^{-\frac{i\tau}{2} {\bT}}\Phi^{[n]})} -I_2)e^{-\frac{i\tau}{2}{\bT}}\Phi^{[n]}\right).
\end{align*}
Similar to the error estimates in \citet{BCJY}, we have the following error bounds
\begin{equation}
\left\|Q^n(x)\right\|_{H^1} \lesssim \eps^2\tau\left(h^{m-1} + \left\|\widetilde{\bf e}^n\right\|_{H^1}\right),	
\end{equation}
Thus, we could obtain
\begin{equation}
\left\|\widetilde{\bf e}^{n+1}\right\|_{H^1} \leq \left\|\widetilde{\bf e}^{n}\right\|_{H^1} +  C_2 \eps^2\tau\left(h^{m-1} + \left\|\widetilde{\bf e}^n\right\|_{H^1}\right), \quad 0 \leq n \leq \frac{T/\eps^2}{\tau}-1,
\end{equation}
where $C_2$ is a constant independent of $h$, $\tau$, $n$ and $\eps$. Since $\widetilde{\bf e}^0 = P_M \Phi_0 - I_M \Phi_0$, we have $\left\|\widetilde{\bf e}^0\right\|_{H^1} \lesssim h^{m-1}$ and discrete Gronwall inequality implies $\left\|\widetilde{\bf e}^{n+1}\right\|_{H^1} \lesssim h^{m-1}$ for $0 \leq n \leq \frac{T/\eps^2}{\tau}-1$. Combining the above estimates with \eqref{eq:diff}, we derive
\[\left\|\Phi(t_n, x) - I_M\Phi^n\right\|_{H^1} \lesssim h^{m-1} + \eps^2 \tau^2+\tau^m_0,\quad 0\leq n\leq \frac{T/\eps^2}{\tau},\]
which shows the improved error bound \eqref{eq:full_bound}. The proof for Theorem \ref{thm:TSFP} is completed.

\section{Numerical results and extensions}
In this section, we present some numerical results of the TSFP method for solving the NLDE \eqref{eq:NLDE_21} in terms of the mesh size $h$, time step $\tau$ and the parameter $0 < \varepsilon \leq 1$ to illustrate our improved uniform error bounds.
\subsection{The long-time dynamics in 1D}  To test the accuracy, we choose the initial data as
\begin{equation}
\phi_1(0, x) = \frac{2}{2+\sin^2(x)},\quad \phi_2(0, x) = \frac{2}{1+\sin^2(x)},\quad x \in (0, 2\pi).
\end{equation}
Since the exact solution is unknown, we use the TSFP method with a very small time step $\tau_e = 10^{-4}$ and a fine mesh size $h_e = \pi/64$ to generate the `reference' solution numerically. Let $\Phi^n = (\Phi^n_0, \Phi^n_1, \ldots, \Phi^n_M)^T$ be the numerical solution obtained by the TSFP method with a given mesh size $h$, time step $\tau$ and the parameter $\eps$ at time $t = t_n$, then we introduce the  discrete $H^1$-error of the wave function as
\begin{align*}
e(t_n) = \|\Phi^n - \Phi(t_n, \cdot)\|_{H^1} 
= \sqrt{h\sum^{M-1}_{j=0} |\Phi^n_j - \Phi(t_n, x_j)|^2+ h\sum^{M-1}_{j=0} |(\Phi')^n_j - \Phi'(t_n, x_j)|^2}.
\end{align*}
where $(\Phi')^n_j$ is defined in \eqref{eq:prime}.

For the long-time dynamics, we quantify the errors as $e_{\max}(t_n) = \max_{0 \leq q \leq n} e(t_q)$. In the rest of the paper, the spatial mesh size is always chosen sufficiently small such that the spatial errors can be neglected when considering the long-time temporal errors.

\begin{figure}[ht!]
\centerline{\includegraphics[width=12cm,height=5cm]{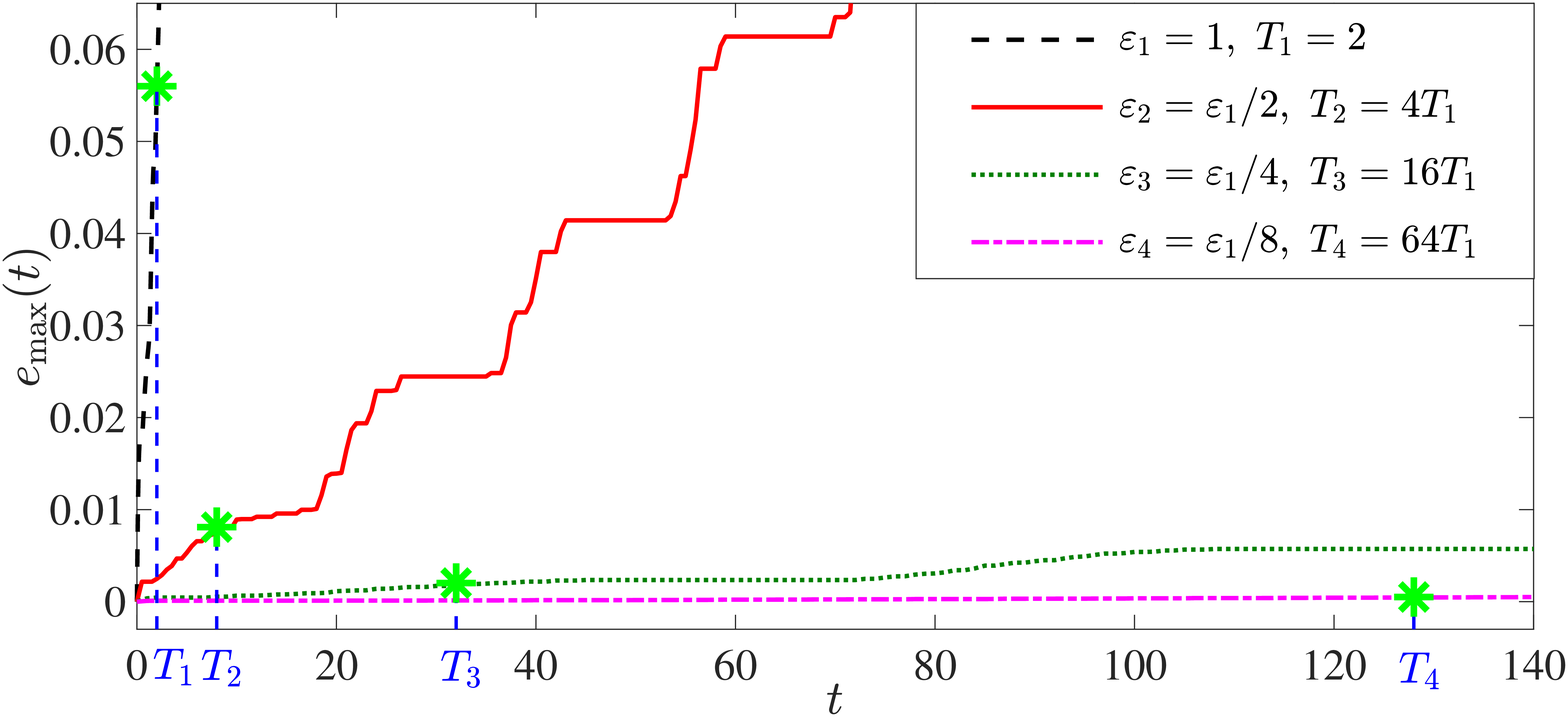}}
\caption{Long-time temporal errors for the wave function of the TSFP method for the NLDE \eqref{eq:NLDE_21} in 1D with different $\eps$.}
\label{fig:1D_long}
\end{figure}

\begin{figure}[ht!]
\begin{minipage}{0.49\textwidth}
\centerline{\includegraphics[width=6cm,height=5cm]{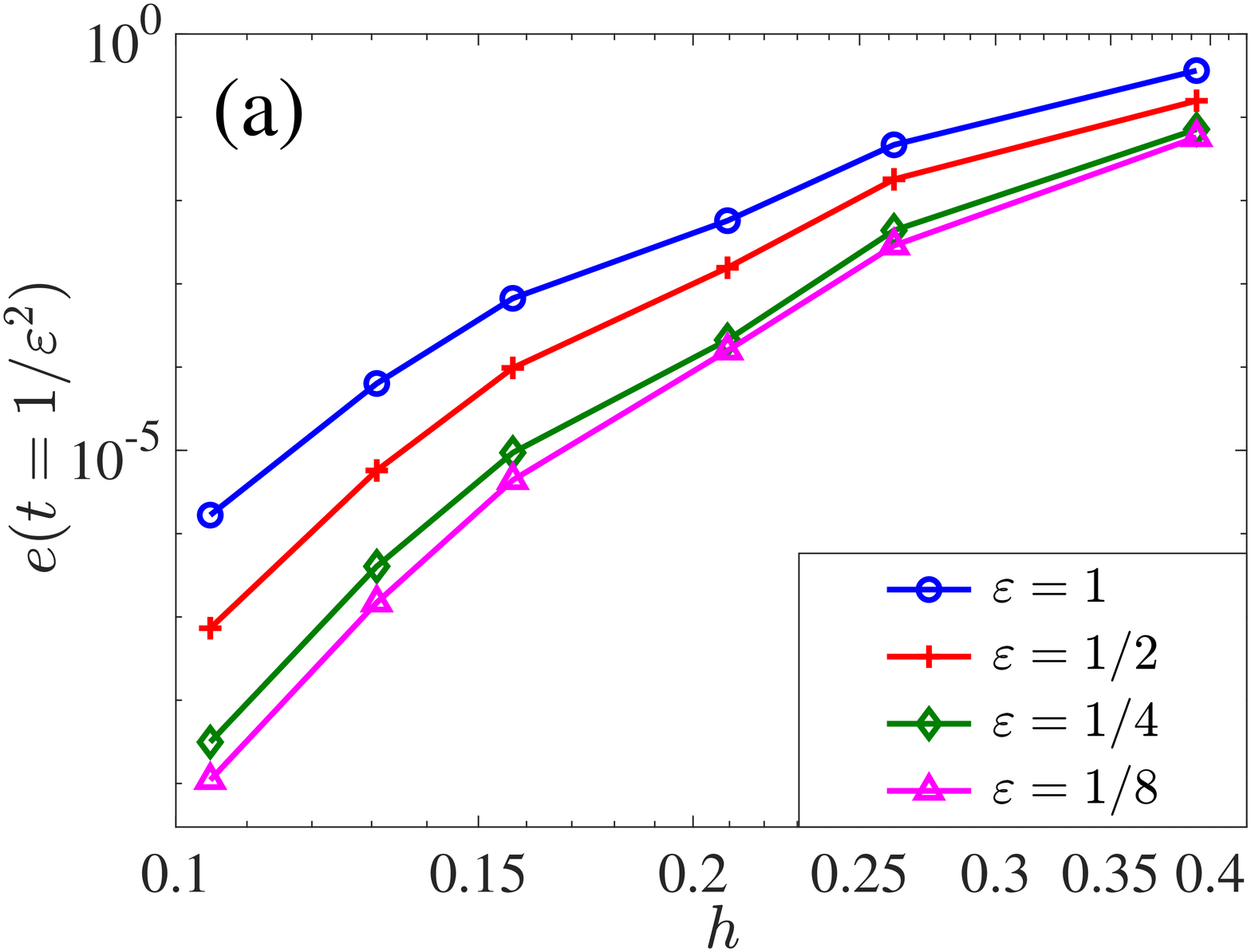}}
\end{minipage}
\begin{minipage}{0.49\textwidth}
\centerline{\includegraphics[width=6cm,height=5cm]{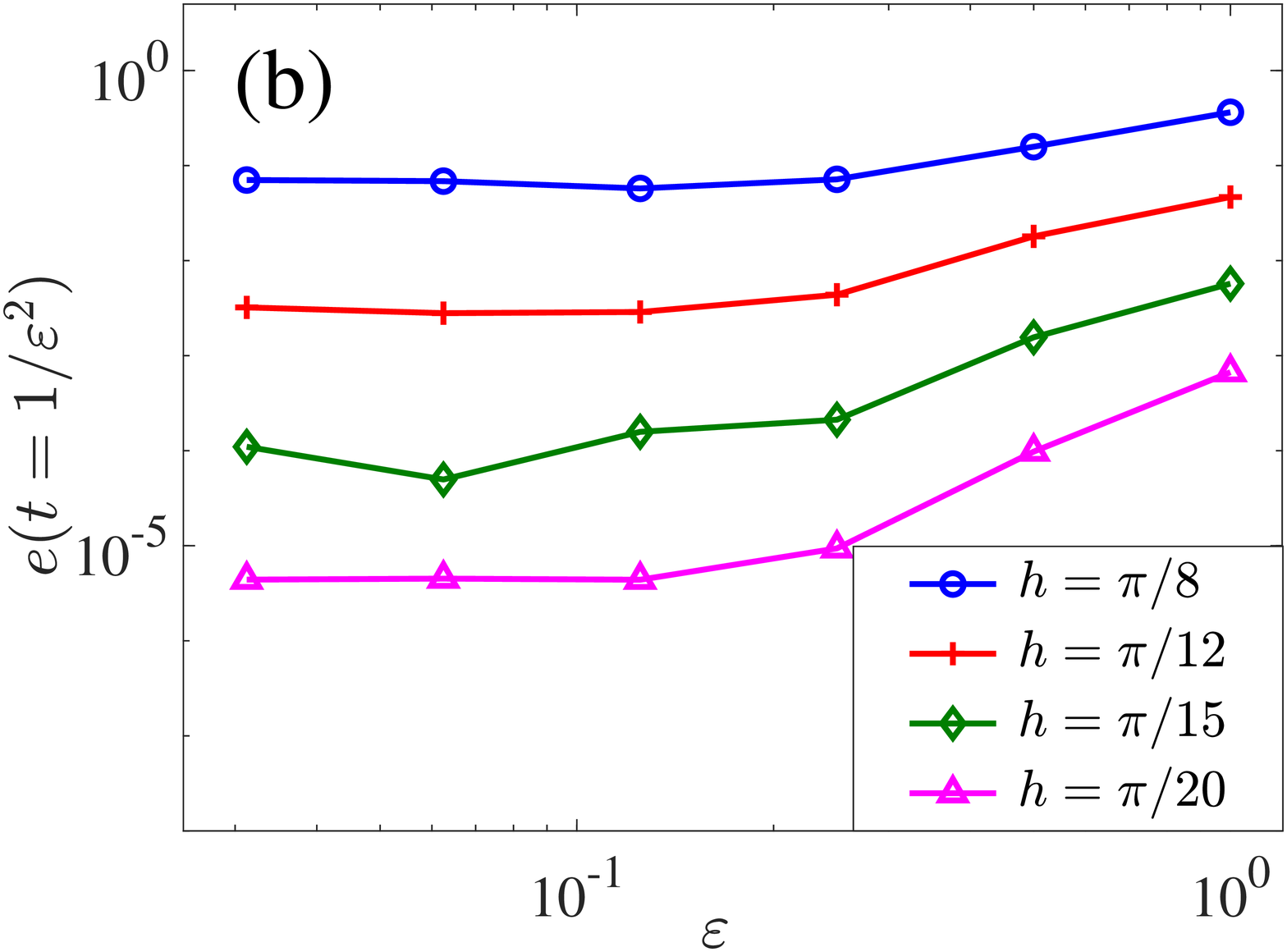}}
\end{minipage}
\caption{Long-time spatial errors of the TSFP method for the NLDE \eqref{eq:NLDE_21} in 1D at $t = 1/\eps^2$.}
\label{fig:1D_spatial}
\end{figure}

\begin{figure}[ht!]
\begin{minipage}{0.49\textwidth}
\centerline{\includegraphics[width=6cm,height=5cm]{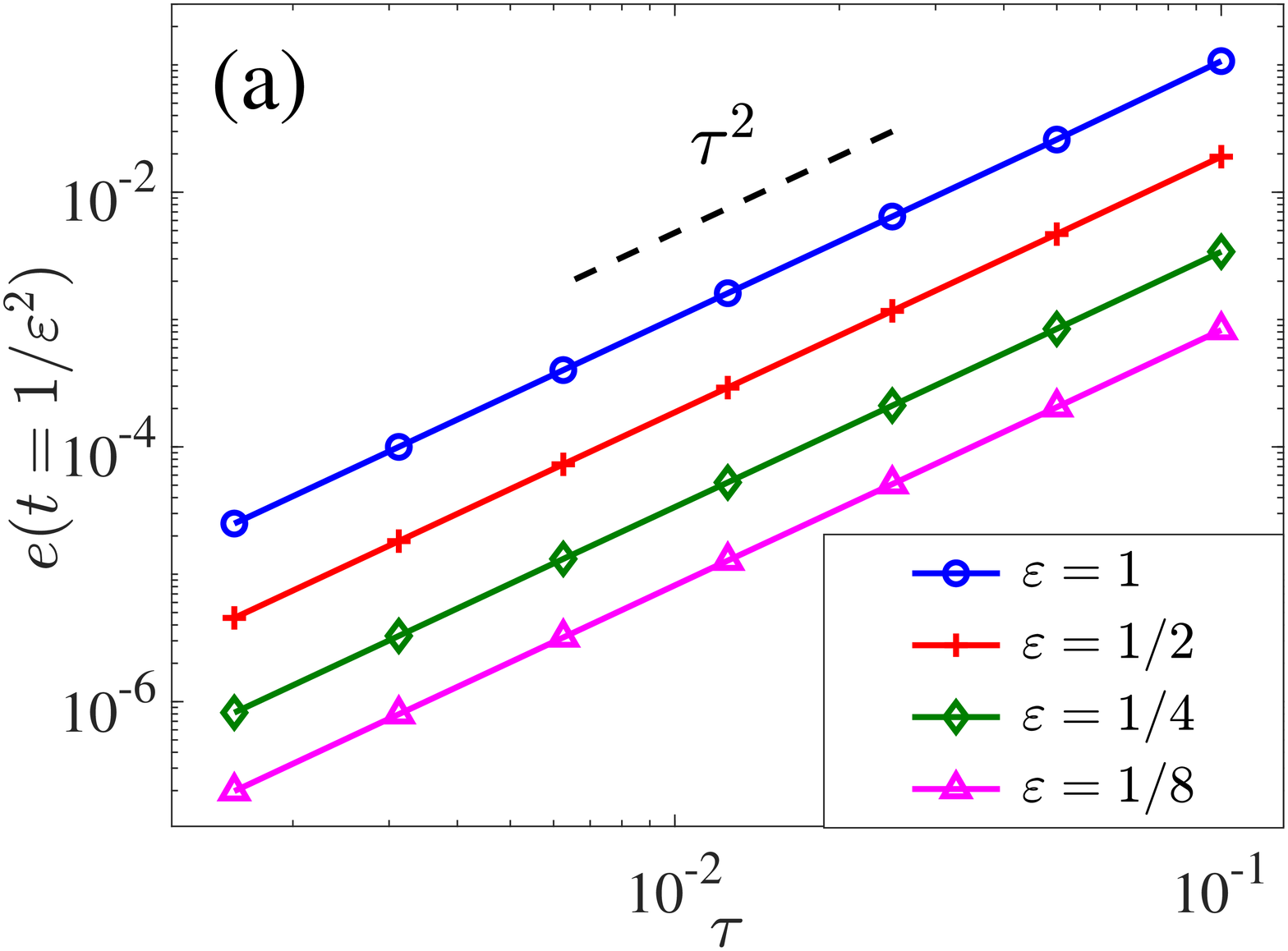}}
\end{minipage}
\begin{minipage}{0.49\textwidth}
\centerline{\includegraphics[width=6cm,height=5cm]{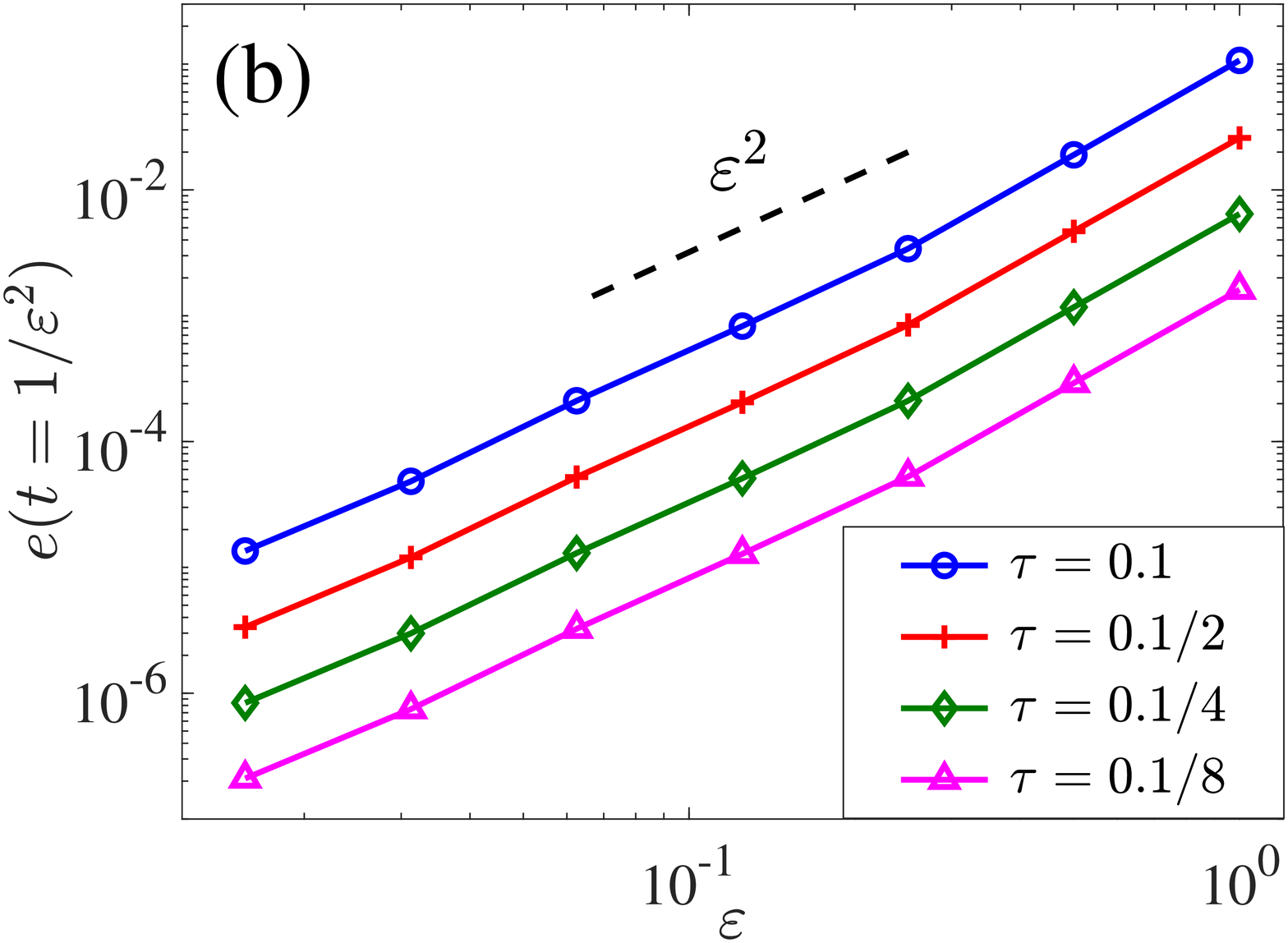}}
\end{minipage}
\caption{Long-time temporal errors of the TSFP method for the NLDE \eqref{eq:NLDE_21} in 1D at $t = 1/\eps^2$.}
\label{fig:1D_temporal}
\end{figure}

\begin{figure}[ht!]
\begin{minipage}{0.5\textwidth}
\centerline{\includegraphics[width=6cm,height=5cm]{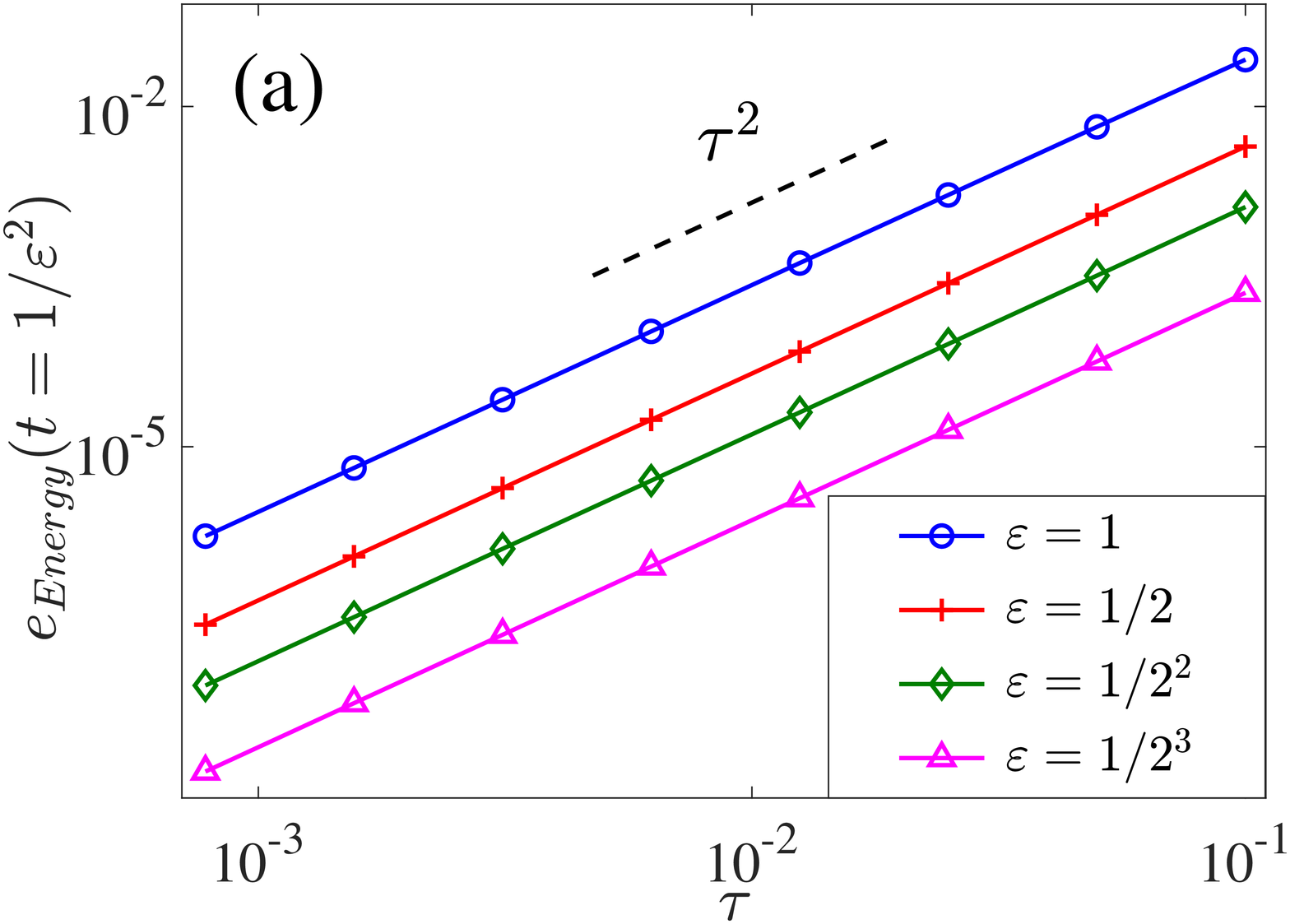}}
\end{minipage}
\begin{minipage}{0.5\textwidth}
\centerline{\includegraphics[width=6cm,height=5cm]{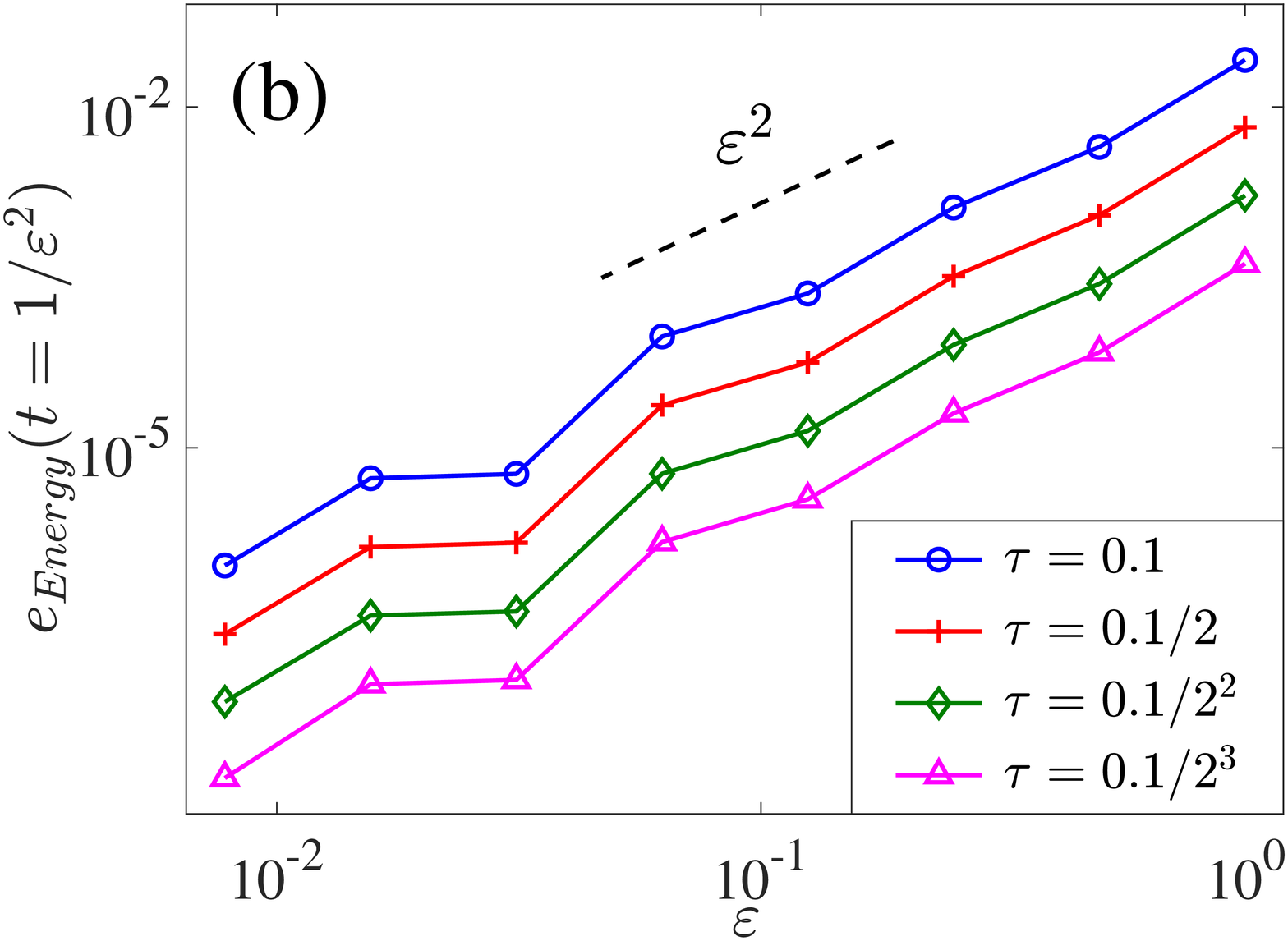}}
\end{minipage}
\caption{Long-time error for the discretized energy $E^n_h$ of the TSFP method for the NLDE \eqref{eq:NLDE_21} in 1D at $t = 1/\eps^2$.}
\label{fig:Energy}
\end{figure}

Fig. \ref{fig:1D_long} displays the long-time errors of the TSFP method \eqref{eq:TSFP} for the NLDE \eqref{eq:NLDE_21} in 1D with different $\eps$, which confirms the improved uniform error bounds at $O(\eps^2\tau^2)$ up to the time at $O(1/\eps^2)$. Fig. \ref{fig:1D_spatial} \& Fig. \ref{fig:1D_temporal} show the spatial and temporal errors of the TSFP methods for the NLDE \eqref{eq:NLDE_21} in 1D at $t = 1/\eps^2$, respectively. Fig. \ref{fig:1D_spatial} indicates the spectral accuracy of the TSFP method in space and the spatial errors are independent of the parameter $\eps$. From Fig. \ref{fig:1D_temporal} (a), we observe the second-order convergence of the TSFP method in time for the fixed $\eps$. Fig. \ref{fig:1D_temporal} (b) again validates that the long-time errors behave like $O(\eps^2\tau^2)$ up to the time at $O(1/\eps^2)$. Fig. \ref{fig:Energy} shows the error of the discrete energy $E^n_h$ also behaves like $O(\eps^2\tau^2)$ up to the time at $O(1/\eps^2)$.

\subsection{The long-time dynamics in 2D} In this subsection, we show an example in 2D with the irrational aspect ratio of the domain $(x, y) \in (0, 2\pi) \times (0, 1)$. In the numerical experiment, we choose the initial data as
\begin{equation}
\phi_1(0, x) = \sin(2x) +\sin(2\pi y), \quad  \phi_2(0, x) = \frac{1}{1 + \cos^2(2x)} +\cos(2\pi y).
\end{equation}

\begin{figure}[ht!]
\centerline{\includegraphics[width=12cm,height=5cm]{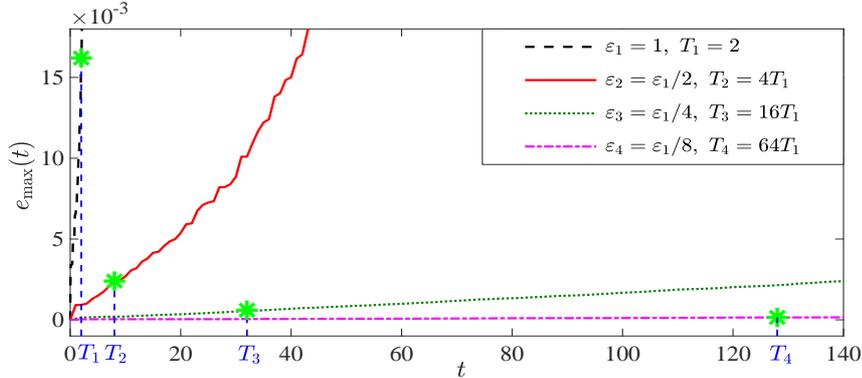}}
\caption{Long-time temporal errors of the TSFP method for the NLDE \eqref{eq:NLDE_21} in 2D with different $\eps$.}
\label{fig:2D_long}
\end{figure}

\begin{figure}[ht!]
\begin{minipage}{0.49\textwidth}
\centerline{\includegraphics[width=6cm,height=5cm]{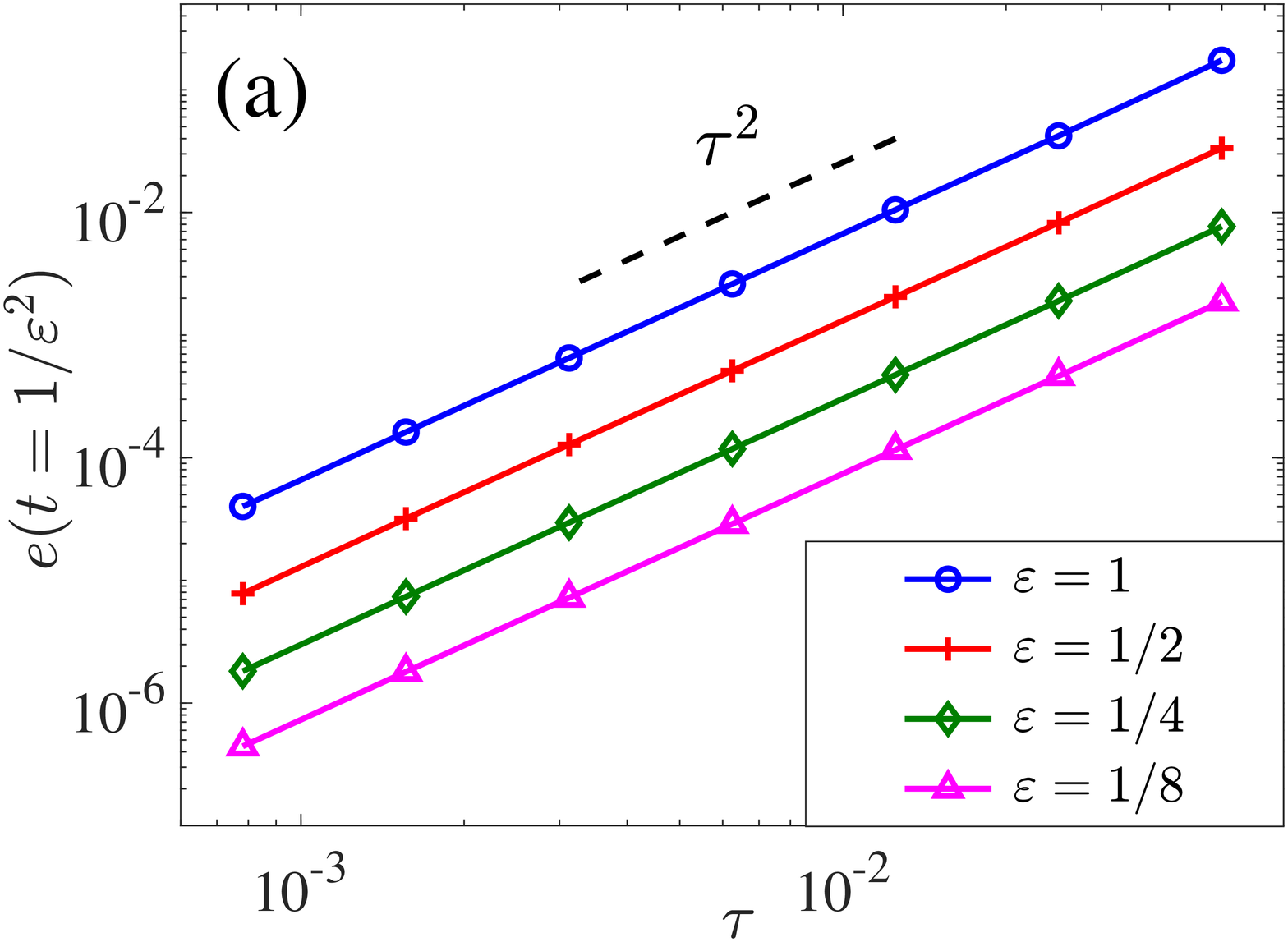}}
\end{minipage}
\begin{minipage}{0.49\textwidth}
\centerline{\includegraphics[width=6cm,height=5cm]{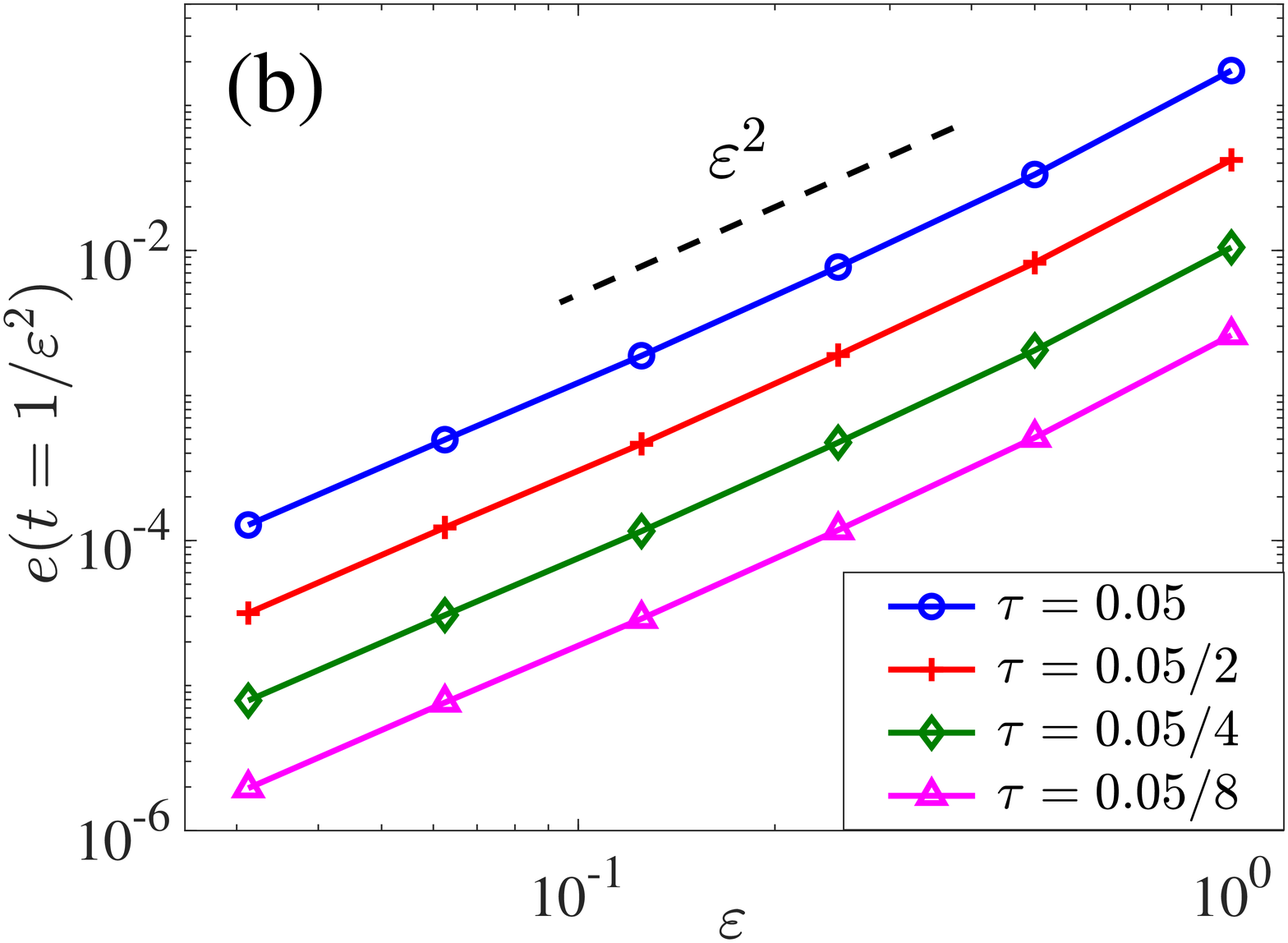}}
\end{minipage}
\caption{Long-time temporal errors of the TSFP method for the NLDE \eqref{eq:NLDE_21} in 2D at $t = 1/\eps^2$.}
\label{fig:2D_temporal}
\end{figure}

Fig. \ref{fig:2D_long} presents the long-time errors of the TSFP method for the NLDE \eqref{eq:NLDE_21} in 2D with a fixed time step $\tau$ and different $\eps$, which confirms the improved uniform error bounds in 2D. Fig. \ref{fig:2D_temporal} plots the temporal errors of the TSFP method for the NLDE in 2D at $t = 1/\eps^2$, which indicates the second-order convergence in time and again validates the improved uniform error bounds at $O(\eps^2\tau^2)$ up to the time at $O(1/\eps^2)$.

\subsection{Extension to the oscillatory NLDE}
In this subsection, we extend the TSFP method and error estimates to an oscillatory NLDE which propagates waves with wave length at $O(\eps^2)$ in time and wave speed at $O(\eps^{-2})$ in space. Introduce a re-scale in time $s = \eps^2 t$ and $\widetilde{\Phi}(s, \bx) = \Phi(t, \bx)$, then the NLDE \eqref{eq:NLDE_21} could be reformulated into the following oscillatory NLDE
\begin{equation}
\label{eq:NLDE_HOE}
i\partial_s\widetilde{\Phi} = \frac{1}{\eps^2}  \Big[- i\sum_{j = 1}^{d}\sigma_j\partial_j + \sigma_3 \Big]\widetilde{\Phi} +  {\bF}(\widetilde{\Phi})\widetilde{\Phi} \quad x \in \Omega, \quad t > 0,
\end{equation}
with the initial data
\begin{equation}
\widetilde{\Phi}(s=0, {\bf x}) = \widetilde{\Phi}_0({\bf x})=O(1), \quad {\bf x} \in \overline{\Omega}.
\end{equation}
The solution of the oscillatory NLDE \eqref{eq:NLDE_HOE} propagates waves with amplitude at $O(1)$, wavelength at $O(1)$ and $O(\eps^2)$ in space and time, respectively, and wave speed at $O(\eps^{-2})$ in space. We remark here the oscillatory nature of the NLDE \eqref{eq:NLDE_HOE} is quite different from the NLDE in the nonrelativistic limit regime, which has been widely studied in \citet{BCJY}, \citet{CW}, \citet{HW}, \citet{KSZ} and \citet{LMZ}. In fact, in this regime, the wave speed is at $O(1/\eps^2)$, while the NLDE in the nonrelativistic limit regime propagates waves with wave speed at $O(1)$. According to the time rescaling, by taking the time step $\kappa = \eps^2 \tau$, we could extend the improved error bounds on the TSFP method for the long-time problem to the oscillatory NLDE \eqref{eq:NLDE_HOE} up to the fixed time $T$. We also just present the result in 1D and it is straightforward to extend to 2D and 3D cases.

\begin{theorem} Let $\widetilde{\Phi}^n$ be the numerical approximation obtained from the TSFP method for the oscillatory NLDE \eqref{eq:NLDE_HOE} in 1D. Assume the exact solution $\widetilde{\Phi}(s, x)$ satisfies
\begin{equation*}
\widetilde{\Phi}(s, x)  \in {L^{\infty}([0, T]; (H^{m}(\Omega))^2)}, \quad m \geq 3,
\end{equation*}
then there exist $h_0 >0$ and $0 < \kappa_0 < 1$ sufficiently small and independent of $\eps$ such that for any $0 < \eps \leq 1$, when $0 < h \leq h_0$ and $0 < \kappa \leq \eps^2\alpha \kappa_0$ for a fixed constant $\alpha \in (0, 1)$, the following improved uniform error bounds hold
\begin{equation}
\left\|\widetilde{\Phi}(s_n, x) - I_M\widetilde{\Phi}^{n}\right\|_{H^1} \lesssim h^{m-1} + \frac{\kappa^2}{\eps^2} + \kappa_0^{m-1},\quad 0 \leq n \leq \frac{T}{\kappa}.
\end{equation}
In particular, if the exact solution is sufficiently smooth, e.g. $\widetilde{\Phi}(t, x) \in (H^{\infty})^2$, the improved uniform error bounds for small $\kappa_0$ would become
\begin{equation}
\left\|\widetilde{\Phi}(s_n, x) - I_M\widetilde{\Phi}^{n}\right\|_{H^1} \lesssim h^{m-1} + \frac{\kappa^2}{\eps^2} , \quad 0 \leq n \leq \frac{T}{\kappa}.
\end{equation}
\label{thm:HOE}
\end{theorem}

The proof of the improved error bounds for the oscillatory NLDE \eqref{eq:NLDE_HOE} in Theorem \ref{thm:HOE} is quite similar to the long-time problem. We omit the details here for brevity and present some numerical results to confirm the sharpness of the improved error bounds. The initial data is chosen as
\begin{equation}
\widetilde{\phi}_1(0, x) = 4x^4(1-x)^4+2,\quad \widetilde{\phi}_2(0, x) = 4x^4(1-x)^4,\quad x \in (0, 1).
\end{equation}
The regularity is enough to ensure the improved error bounds.

\begin{table}[h!]
\caption{Temporal errors of the TSFP method for the NLDE \eqref{eq:NLDE_HOE} in 1D at $t = 1$.}
\label{tab:HOE}
\renewcommand{\arraystretch}{1.3}
\centering
\vspace{0.3cm}
\begin{tabular}{cccccc}
\hline
$e(t=1)$ &$\kappa_0 = 0.05$ & $\kappa_0/4 $ &$\kappa_0/4^2 $ & $\kappa_0/4^3$ & $\kappa_0/4^4$  \\
\hline
$\varepsilon_0 = 1$ & \bf{1.26E-2} & 7.57E-4 & 4.72E-5 & 2.95E-6 & 1.84E-7 \\
order & \bf{-} & 2.03 & 2.00 & 2.00 & 2.00  \\
\hline
$\varepsilon_0 / 2$ & 1.21E-1 & \bf{5.71E-3} & 3.52E-4 & 2.20E-5 & 1.37E-6  \\
order & - & \bf{2.20} & 2.01 & 2.00 & 2.00 \\
\hline
$\varepsilon_0 / 2^2$ & 5.39E-2 & 9.18E-3 & \bf{3.85E-4} & 2.37E-5 & 1.48E-6  \\
order & - & 1.23 & \bf{2.29} & 2.01 & 2.00  \\
\hline
$\varepsilon_0 / 2^3$ & 2.40E-1 & 1.08E-2 & 2.13E-3 & \bf{8.68E-5} & 5.34E-6  \\
order & - & 2.24 & 1.17 & \bf{2.31} & 2.01 \\
\hline
$\varepsilon_0 / 2^4$ & 1.17E-1 & 3.80E-2 & 2.46E-3 & 6.76E-4 & \bf{2.46E-5}  \\
order & - & 0.81 & 1.97 & 0.93 & \bf{2.39}  \\
\hline
\end{tabular}
\end{table}

Table \ref{tab:HOE} lists the temporal errors of the TSFP method for the oscillatory NLDE \eqref{eq:NLDE_HOE} in 1D with different $\eps$. It can be clearly observed that the second-order convergence can only be observed when $\kappa \lesssim \eps^2$ (cf. the upper triangle above the diagonal with bold letters) and the temporal errors behave like $O(\kappa^2/\eps^2)$. Along each diagonal in the table, i.e. with $\kappa=O(\eps^2)$, we observe linear convergence with respect to $\kappa$, which confirm the error bound at
$\frac{\kappa^2}{\eps^2}=O\left(\frac{\kappa^2}{\kappa}\right)=O(\kappa)$ 
in \eqref{thm:HOE}.

%%%%%%%%%%%%%%%%%%%%%%%%%%%%%%%%%%
% %    Section 6 Conclusions
%%%%%%%%%%%%%%%%%%%%%%%%%%%%%%%%%%
\section{Conclusions}
Improved uniform error bounds on time-splitting methods for the long-time dynamics of the weakly nonlinear Dirac equation (NLDE) were rigorously proven. With the help of the regularity compensation oscillation (RCO) technique, the long-time errors up to the time at $O(1/\eps^2)$ for the semi-discretization and full-discretization are at $O(\eps^2\tau^2)$ and $O(h^{m-1} + \eps^2\tau^2)$, respectively, which improve the standard error bounds in the literature at $O(\tau^2)$ and $O(h^{m-1} +\tau^2)$ for the semi-discretization and full-discretization, respectively, especially when $0<\eps\ll1$. The improved error bounds were extended to an oscillatory NLDE up to a fixed time $T$. Numerical results in 1D and 2D agreed well with the theoretical results and suggested that the error bounds are sharp.

\section*{Acknowledgements}
This work was partially supported by the Ministry of Education of Singapore grant MOE2019-T2-1-063 (R-146-000-296-112, W. Bao, Y. Feng) and Natural Science Foundation of China grant 12171041 and 11771036 (Y. Cai).

\end{document}